\documentclass{amsart}
\usepackage[backref,colorlinks=true,
            linkcolor=blue,
            citecolor=blue,
            urlcolor=blue]{hyperref}
\usepackage{amssymb,amsmath}
\usepackage{graphicx}%
\graphicspath{{../img_pdf/}}

\usepackage{multirow}%
\usepackage{amsmath,amssymb,amsfonts}%
\usepackage{amsthm}%
\usepackage{mathrsfs}%
\usepackage[title]{appendix}%
\usepackage{xcolor}%
\usepackage{textcomp}%
\usepackage{manyfoot}%
\usepackage{booktabs}%
\usepackage{algorithm}%
\usepackage{algorithmicx}%
\usepackage{algpseudocode}%
\usepackage{listings}%
\usepackage{overpic}%

\usepackage{bm,color}
\usepackage{pifont}
\usepackage{enumitem}
\usepackage{subcaption}

\def\Z{\mathbb{Z} } 
 
\def\R{\mathbb{R} }

\newcommand{\B}{\Box}

\theoremstyle{definition}
\newtheorem{definition}{Definition}
\newtheorem{theorem}{Theorem}
\newtheorem{proposition}[theorem]{Proposition}
\newtheorem{lemma}[theorem]{Lemma}

\newtheorem{remark}{Remark}

\title[Topological vortex identification for 2D turbulent flows]{Topological vortex identification for two-dimensional turbulent flows in doubly periodic domains}

\author[M. Kimura]{Mitsuaki Kimura}
\address{Department of Mathematics, Osaka Dental University, 
8-1 Kuzuha Hanazonocho, Hirakata, Osaka 573-1121, Japan}

\author[T. Matsumoto]{Takeshi Matsumoto}
\address{Department of Physics, Kyoto University, 
Kitashirakawa Oiwake-cho, Sakyo-ku, Kyoto 606-8502, Japan}

\author[T. Sakajo, H. Takeuchi]{Takashi Sakajo \and Hiroshi Takeuchi}
\address{Department of Mathematics, Kyoto University, 
Kitashirakawa Oiwake-cho, Sakyo-ku, Kyoto 606-8502, Japan}

\author[T. Yokoyama]{Tomoo Yokoyama}
\address{Department of Mathematics, Saitama University, 
255 Shimo-Okubo, Sakura-ku, Saitama 338-8570, Japan}

\date{\today}


\begin{document}

\maketitle

\begin{abstract}
The dynamics and statistical properties of two-dimensional (2D) turbulence are often investigated through numerical simulations of incompressible, viscous fluids in doubly periodic domains. 
A key challenge in 2D turbulence research is accurately identifying and describing statistical properties of its coherent vortex structures within complex flow patterns. 
This paper addresses this challenge by providing a classification theory for the topological structure of particle orbits generated by instantaneous Hamiltonian flows on the torus $\mathbb{T}^2$, which serves as a mathematical model for 2D incompressible flows.
Based on this theory, we show that the global orbit structure of any Hamiltonian flow can be converted into a planar tree, named a partially Cyclically-Ordered rooted Tree (COT), and its corresponding string expression (COT representation). 
We apply this conversion algorithm to 2D energy and enstrophy cascade turbulence. 
The results show that the complex topological structure of turbulent flow patterns can be effectively represented by simple trees and sequences of letters, thereby successfully extracting coherent vortex structures and investigating their statistical properties from a topological perspective.
\end{abstract}

\section{Introduction}\label{sec1}
We consider the flow generated by an inviscid and incompressible fluid of constant density in the flat torus $\mathbb{T}^2$, that is, a two-dimensional domain with a doubly periodic boundary condition.
The evolution of the velocity field ${\bm u}({\bm x}, t)=(u(x,y,t), v(x,y,t))$ at ${\bm x} = (x,y) \in \mathbb{T}^2$ and time $t \in \mathbb{R}$ is governed by the Navier-Stokes equation.
\[
\partial_t {\bm u} + ({\bm u} \cdot \nabla) {\bm u} = -\nabla p + \nu \triangle {\bm u} + {\bm f}, \qquad \nabla \cdot {\bm u} = 0,
\]
in which $p({\bm x}, t)$ and ${\bm f}$ denote the pressure and an external force. 
For a smooth velocity field, we define the stream function $\psi({\bm x},t)$ by $(u,v)=(\partial_y \psi, -\partial_x \psi)$ and the vorticity $\omega({\bm x}, t)= \partial_x v -\partial_y u$, respectively; the Navier--Stokes equation is reduced to 
\[
\partial_t \omega + \mathcal{J}(\omega, \psi)=  \nu \triangle \omega + \widetilde{f}, \qquad -\triangle \psi = \omega,
\]
where $\mathcal{J}(f,g)$ denotes the Jacobi matrix for the functions $f$ and $g$, and a scalar function $\widetilde{f}$ is defined by $\nabla \times {\bm f}=(0,0,\widetilde{f})$.

These fluid equations on the flat torus are commonly used to investigate two-dimensional turbulent flows numerically to clarify the relation between their statistical properties and complex flow structures.
This is because accurate schemes such as the spectral method are available. 
See a comprehensive survey paper by Tabeling~\cite{Tabeling2002-by} on the physics of two-dimensional turbulence. 
Typically, two types of two-dimensional turbulence are studied. 
One is called free decaying turbulence. 
It describes how a complex flow relaxes into a simpler flow pattern without external force $\widetilde{f}=0$ while losing its energy due to viscous dissipation. 
The other is known as inverse cascade turbulence. 
When an external force is applied to a flow, we investigate the ensemble average of the energy spectra, which is the $L^2$ norm of the velocity field, and the enstrophy spectra, which is the $L^2$ norm of vorticity. 
We also observe that the complex flow pattern with many coherent vortex structures gradually transitions to that with larger-scale vortices. 
The decaying process is characteristic of two-dimensional turbulence and is often investigated as the interaction among coherent vortex structures of various scales in complex turbulent flow patterns.
Therefore, extracting coherent vortex structures from complex flow patterns is an essential task in the study of two-dimensional turbulence.
  
Various methods have been proposed to extract the coherent rotational component from the vorticity distribution. 
For example,  McWilliams~\cite{Mcwilliams1990} proposed an intuitive method of finding the point where the vorticity is at its extreme value and extracting a closed contour enclosing the vorticity region around the point.
Basdevant~\cite{Basdevant1994} and Hua \& Klein~\cite{Hua1998} defined a rotational flow component using the eigenvalues of the velocity gradient tensor based on the Weiss assumption that it varies very slowly.
Farge~\cite{Farge1992} used the wavelet transform of the vorticity and coherent rotating regions, which have been extracted corresponding to the coefficients of an appropriate scale. 
All of these methods extract coherent vortex structures from the vorticity function.
However, since vorticity has large values not only in rotating flow regions but also in strong shear regions, it is still unclear how we define a vortex region. 
Hence, it is inevitable to introduce certain empirical hyperparameters and working hypotheses for vortex identifications in these methods.

As an alternative approach to extracting coherent rotating flow structures from two-dimensional flow patterns, a method called Topological Flow Data Analysis (TFDA) has been developed based on topology and the theory of dynamical systems. 
Recall that the incompressible flow field constructed here is a Hamiltonian vector flow with the steam function $\psi$ being its Hamiltonian.
Then, the particle orbits (streamlines) generated by a stationary Hamiltonian vector field coincide with the level curves of the Hamiltonian due to $\nabla \psi \cdot {\bm u}=0$. 
Therefore, a topological classification of the Hamiltonian flow is equivalent to a classification of the topological structure of the level curves. 
Ma and Wang~\cite{ma2005geometric} considered such a time-independent Hamiltonian vector field on a compact surface. 
They show that the topological orbit structures consist of a finite number of topologically distinguishable local flow patterns under the condition of structural stability in $C^r$ ($r\geqq 1$) topology.
In~\cite{yokoyama2013word}, it is shown that this topological classification theory can be extended to Hamiltonian flows in the presence of uniform flows.
Later, it was discovered that the global topological structure is converted uniquely into a planar tree~\cite{sakajo2018tree}, which is now called a partially cyclically ordered rooted tree (COT), and its symbolic representation is called the COT representation~\cite{yokoyama2021cot}. 
The COT is found to be isomorphic to the Reeb graph of Hamiltonian as an abstract graph~\cite{uda2019persistent_en}. 
Consequently, TFDA allows us to distinguish the topological structures of flow patterns with COT and COT representation.
Using the COT, we can extract the area surrounded by a self-connected saddle separatrix as a coherent rotational flow component.
For example, by tracking the evolution of coherent cyclonic rotating flow regions extracted using this method from atmospheric data on 500-hPa height surfaces, we identify the occurrence of atmospheric blocking phenomena~\cite{UDA2021}.
In oceanography, it is known that the Kuroshio, a flow along the southern coast of Japan, exhibits a meandering phenomenon called the Kuroshio Large Meander (KLM). 
By extracting the coherent eddy structure due to the meandering from the sea surface height data with TFDA, we have successfully estimated the duration of KLM, which agrees well with the empirical prediction owing to the Japan Meteorological Agency.
Note that this topological classification theory using COT and COT representations is extended to that for a wide class of vector fields, including compressible flows~\cite{Sakajo2022-DDMA}, and it is used effectively in the identification of vortex structures for blood flows in the left ventricle of the heart~\cite{Sakajo2023-SIAM}. 

In this paper, we extend this topological flow data analysis to Hamiltonian vector fields on $\mathbb{T}^2$, so that it is applicable to extracting coherent vortex structures from complex flow patterns in two-dimensional turbulence.
Since the flat torus is compact and its topology is different from $\mathbb{R}^2$, the original TFDA theory for Hamiltonian flows on $\mathbb{R}^2$ cannot be applied directly as is.
Since an incompressible vector field on a compact surface can be described locally in terms of differential forms in general, it is impossible to provide a unique Hamiltonian and identify its level curves for the surface.
In other words, by using the properties of a manifold, a Hamiltonian is locally constructed on a subset of $\mathbb{R}^2$ as its image for each local coordinate system, and it is impossible to select a single Hamiltonian.
However, in the case of $\mathbb{T}^2$, all points can be uniquely represented in a single local coordinate system if we identify this local coordinate system with the $(x,y)$-plane.
Hence,  we define a Hamiltonian $H(x,y)$ in that plane under the assumption that the harmonic part is absent. 
Hence, the orbits coincide with the level curves of this Hamiltonian similar to $\mathbb{R}^2$.
Another mathematical issue arises when we develop the TFDA theory for the flat torus. 
 In the TFDA theory in $\mathbb{R}^2$, the COT is isomorphic to the Reeb graph of the Hamiltonian. 
However, in general, the Reeb graph of a function on the torus has one loop, which is not a tree (an acyclic graph). 
Hence a device is required to convert the topological structures of the Hamiltonian flows on $\mathbb{T}^2$ into a COT and its associated COT representation.

The purpose of this paper is to establish a TFDA method for flows on a flat torus by addressing these mathematical difficulties.
Then, we apply this method to extract vortex structures in two-dimensional turbulence, thus clarifying their statistical properties.
This paper is structured as follows. 
Section~\ref{sec:2} provides some mathematical preliminaries required to describe our new TFDA theory.
In Section~\ref{sec:3}, we present a topological classification theory for Hamiltonian flows on a torus without a physical boundary.
Based on classification theory, in Section~\ref{sec:4}, we describe a conversion algorithm to COT and COT representations for Hamiltonian flows on a torus. 
Section~\ref{sec:5} shows the applications of our TFDA method to two-dimensional turbulent flows; 
In Section~\ref{sec:5.1}, we demonstrate how our TFDA algorithm can be applied to snapshots of free-decaying turbulence and enstrophy cascade turbulence.
In Section~\ref{sec:5.2}, we apply TFDA to the time evolution of energy cascade turbulence and enstrophy cascade turbulence. 
We elucidate the differences between these two turbulent states in terms of the statistical properties of {\it terminal vortices}, which are defined by COT from a topological perspective.
Section~\ref{sec:6} is a summary. 

\section{Mathematical preliminaries}~\label{sec:2}
We consider a two-dimensional compact manifold $S$.
The flow $v$ on the surface $S$ is a continuous $\mathbb{R}$-action on $S$, that is, $v \colon \mathbb{R} \times S \rightarrow S$. 
For the flow $v$, we define a map $v_s \colon S \rightarrow S$,  where $v_s(x) := v(s, x)$ for $s \in \mathbb{R}$. 
For a point $x \in S$, 
\[
O(x)=\{ v_s(x) \in S \, \vert \, s \in \mathbb{R}\}
\]
is called an orbit passing through $x$.  
We classify the topological structure of the set of all orbits, $\{ O(x) \vert x \in S\}$.
The flow $v$ is an abstract mathematical object, but is regarded as particle orbits generated by a two-dimensional steady incompressible velocity field in fluid mechanics. 
This is because the incompressible vector field is a Hamiltonian vector filed with the stream function being its Hamiltonian, and particle orbits are equivalent to streamlines, i.e., level-curves of the stream function, in this case.
On the other hand, since the stream function (Hamiltonian) is generally time dependent in complex turbulent flows,  these particle orbits do not coincide with the orbit $O(x)$ defined here. 
However, in TFDA, we characterize the evolution of time-dependent Hamiltonian flows as the topological change of Hamiltonian flow patterns for every fixed time.
Hence, in this paper, we regard the orbits in $O(x)$ equivalently as streamlines of instantaneous Hamiltonian, as long as there is no confusion.
 
We assume that all orbits are {\it proper}, which is a notion of the regularity of orbits in foliation theory~\cite{hector1983introduction}.
Although its mathematical definition is not specified here, it is not a substantial limitation for data analysis since this assumption is satisfied by flow data with finite resolution obtained by numerical simulations or measurements. 
It is known that proper orbits are further classified into the following categories.
\begin{definition}
For a proper orbit of the flow $v$ on a surface $S$, singular orbits, periodic orbits, and non-closed orbits are defined as follows.
\begin{itemize}
\item $x\in S$ is a singular orbit if $x=v_s(x)$ for all $s \in \mathbb{R}$, that is, $O(x)=\{ x \}$.
\item An orbit $O(x)$ is periodic if there exist $T>0$ such that $v_T(x)=x$ and $v_s(x) \neq x$ for $0 < s < T$.
\item A proper orbit that is neither a singular point nor a periodic orbit is called a non-closed orbit.
\end{itemize}
\end{definition}
Recall some definitions of special orbits in the flow on $S$.
\begin{definition}
For a proper orbit through $x \in S$ of the flow $v$ on a surface $S$, the $\omega$-limit set $\omega(x)$ and the $\alpha$-limit set $\alpha(x)$ are defined as follows.
\begin{itemize}
\item $\omega(x):=\cap_{n\in \mathbb{R}}\overline{ \{ v_s(x) \, \vert \, s> n\}}$; (The set that $O(x)$ can reach in $s\rightarrow \infty$),
\item $\alpha(x):= \cap_{n \in \mathbb{R}}\overline{ \{ v_s(x) \, \vert \, s<n \}}$; (The set that $O(x)$ can reach in $s\rightarrow -\infty$),
\end{itemize}
where $\overline{A}$ denotes the closure of the subset $A\subseteq S$. 
\end{definition}

Let $I \subset \mathbb{R}$ be an interval and consider a continuous map $\gamma \colon I \to S$ on a surface $S$.  
Then the image $\gamma(I)$ is called a curve on $S$.
Moreover, a curve $\gamma \colon [a,b] \to S$ is said to be simple closed if $\gamma(a)=\gamma(b)$ and $\gamma(s) \neq \gamma(t)$ for every $s,t \in [a,b)$ with $s\neq t$.
The orbit of the flow $v$ in $S$ defines a curve in $S$. 
In particular, a periodic orbit $O(x)$ with period $T$ is a simple closed curve $\gamma:[0,T]\to S$ due to $\gamma(s)=v_s(x)$. 

We now introduce the concepts of \textit{essential} and \textit{inessential} for subsets $S \subseteq \mathbb{T}^2$ on a torus, which plays an important role in the description of the TFDA theory on the torus. 
\begin{definition} \label{def:iness_arb}
Let $S \subseteq \mathbb{T}^2$ be a compact surface. A subset $E \subseteq S$ is called {\it inessential} if there exists an open neighborhood $U$ of $E$ such that every simple closed curve in $U$ is contractible in $\mathbb{T}^2$.
In addition, a subset $E \subseteq S$ is called {\it quasi-inessential} if $\overline{E}$ is inessential.
\end{definition}
An essential set is defined as the negation of an inessential set.
\begin{definition} \label{def:ess_arb}
Let $S \subseteq \mathbb{T}^2$ be a compact surface. A subset $E \subseteq S$ is called {\it essential} if every open neighborhood $U$ of $E$ contains a simple closed curve in $U$ that is not contractible in $\mathbb{T}^2$. In addition, a subset $E \subseteq S$ is called {\it quasi-essential} if $\overline{E}$ is essential.
\end{definition}
Furthermore, we also introduce a subset of essential sets as follows.
\begin{definition} \label{def:fully_ess_arb}
Let $S \subseteq \mathbb{T}^2$ be a compact surface. A subset $E \subseteq S$ is called {\it fully essential} if its complement $\mathbb{T}^2 - E$ is inessential.
\end{definition}
We can also define sets that are essential but not fully essential.
\begin{definition} \label{def:singly_ess_arb}
Let $S \subseteq \mathbb{T}^2$ be a compact surface. When a subset $E \subseteq S$ is essential but not fully essential, it is called {\it singly essential}.
\end{definition}

\begin{figure}[t]
\begin{center}
\includegraphics[bb=0 0 570 455, scale=0.6]{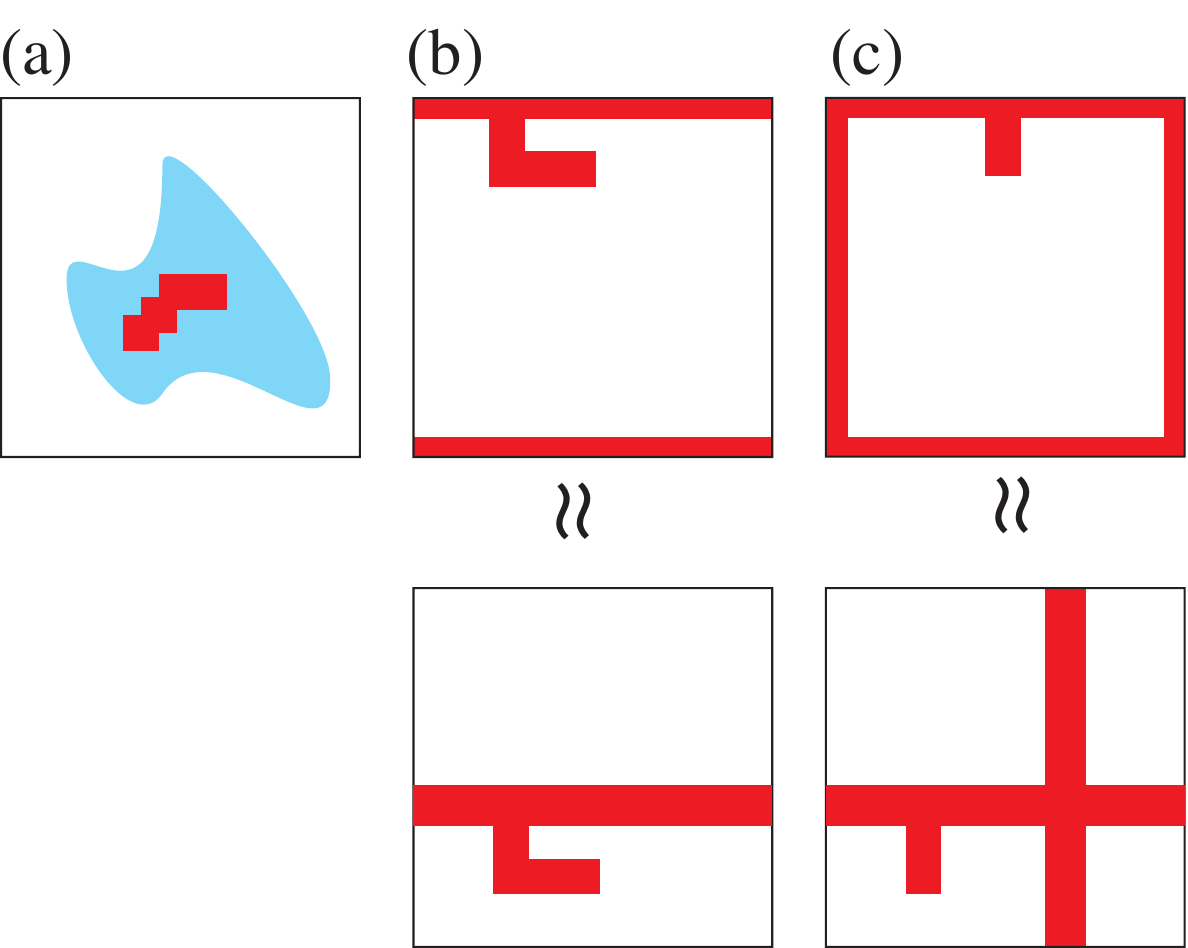}
\end{center} 
\caption{The red filled regions $U$ are (a) essential, (b) singly essential, and (c) fully essential, respectively. Note that the panels in (b) and (c), both in the upper and lower rows, represent the same region under the assumption of doubly periodic boundary condition.}
\label{fig:ess_iness}
\end{figure}

For intuitive understanding, some examples are shown in Figure~\ref{fig:ess_iness}. 
The red filled region $U$ in Figure~\ref{fig:ess_iness}(a) is inessential as it is contained in the disk (or its homeomorphic shape).  
The red subset in Figure~\ref{fig:ess_iness}(b) is singly essential, and its boundary contains an essential simple closed curve. 
The fully essential subset in Figure~\ref{fig:ess_iness}(c) crosses both sides of the periodic boundary of the torus. Each component of that boundary is an inessential simple closed curve. 

\begin{remark}
The terms ``essential'', ``inessential'' and ``fully essential'' are well known in topology. 
On the other hand, ``quasi-essential'' and ``quasi-inessential'' are new terms required to describe topological flow properties. 
Also, ``singly essential'' is a new term, since simple closed curves that are ``essential'' but not ``fully essential'' appear when we consider sets in the torus $\mathbb{T}^2$ in particular. 
Hence, we have the following decomposition for a compact surface $S \subseteq \mathbb{T}^2$:
\[
\begin{split}
2^S 
 = \{ \text{inessential subset of }S \} &\sqcup \{ \text{singly essential subset of }S \} \\
 &\sqcup \{ \text{fully essential subset of }S \}.
\end{split}
\]
Here, $2^S$ denotes the family of sets consisting of all subsets of $S$.
\end{remark}

\section{Topological classification of streamlines in Hamiltonian flows}\label{sec:3}
\subsection{Local orbit structures}~\label{sec:3.1}
A flow generated by a Hamiltonian vector field is called a Hamiltonian flow.
We introduce local orbit structures that are to be topologically identified in Hamiltonian flows.
As explained in Introduction, the Hamiltonian on a torus is represented by a single function, $H(x,y)$, where the local coordinate is identified with the $(x,y)$-plane. 
Therefore, it is sufficient to consider the local topological structure of the level curves of the function $H(x,y)$.

\begin{figure}[t]
\begin{center}
\includegraphics[scale=0.6]{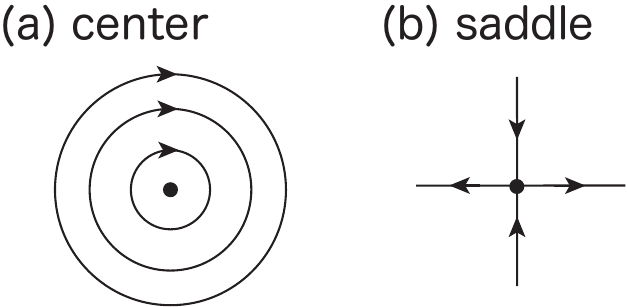}
\end{center}
\caption{Zero-dimensional structures. (a) A center. (b) A saddle.}
\label{fig:Point_Flow_Components}
\end{figure}

Figure~\ref{fig:Point_Flow_Components} shows non-degenerate zero-dimensional singular orbits of a Hamiltonian flow on a torus.
A center is a singular orbit associated with infinitely many periodic orbits around it, as shown in Figure~\ref{fig:Point_Flow_Components}(a).
Figure~\ref{fig:Point_Flow_Components}(b) is a saddle. 
It is associated with four separatrices that are $\alpha$-limit and $\omega$-limit sets of the saddle. 
The separatrix is then said to be connected. 

\begin{figure}[t]
\begin{center}
\includegraphics[bb=0 0 836 204,scale=0.427]{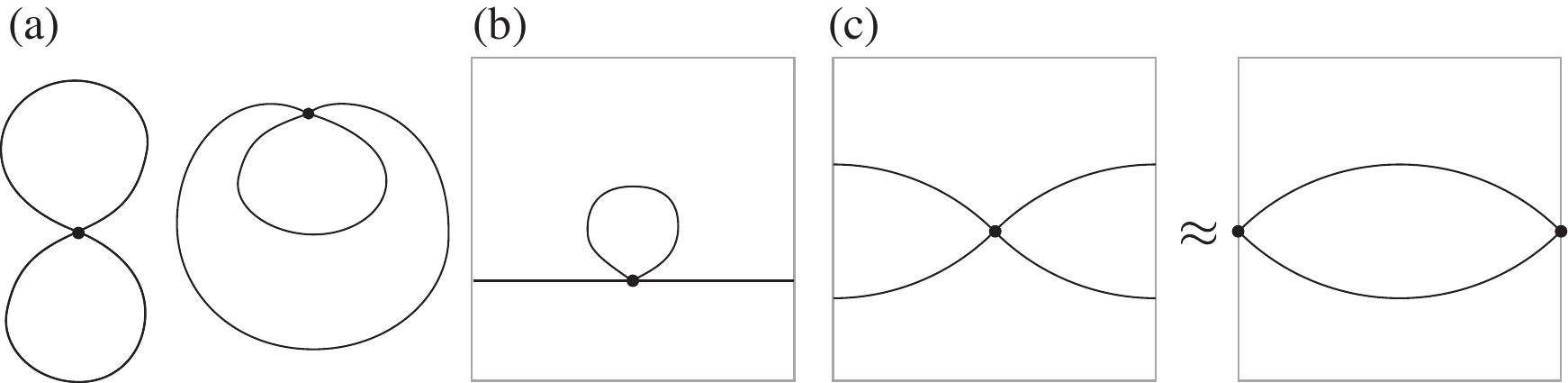}
\end{center} 
\caption{Non-trivial and non-periodic circuits. (a) Saddle connections consisting of two quasi-inessential self-connected saddle separatrices. (b) Saddle connection consisting of a quasi-inessential and a quasi-essential saddle connections. (c) Saddle connection consisting of two quasi-essential saddle separatrices.}
\label{fig:ess_saddle_connection}
\end{figure}

When a saddle separatrix connects the same saddle, it is called self-connected. 
Figure~\ref{fig:ess_saddle_connection} shows examples of self-connected saddle separatrices. 
{\it A saddle connection diagram} for a Hamiltonian flow is the set of saddles and the saddle separatrices that connect them.
In particular, a connected component of a saddle connection diagram is called a {\it saddle connection}.
The saddle connection diagram is said to be {\it self-connected} if it consists of self-connected saddle connections. 

A {\it circuit} is an immersion of a circle in $S \subset \mathbb{T}^2$.
The image of the circle under this mapping is either a point or an orbit homeomorphic to the immersed circle.
When the image is a point, it is called a trivial circuit.  
It is called a non-trivial circuit when it is an immersed circle.
Furthermore, non-trivial circuits are classified into two types. 
One is a periodic orbit, called a {\it cycle}.
The other is a non-trivial circuit consisting of saddles and separatrices.
Figure~\ref{fig:ess_saddle_connection}(a) shows circuits that consist of a saddle and two quasi-inessential saddle separatrices.
Figure~\ref{fig:ess_saddle_connection}(b) and (c) show circuits with singly essential saddle connections.
The orbit structure in Figure~\ref{fig:ess_saddle_connection}(b) (resp. Figure~\ref{fig:ess_saddle_connection}(c)) consists of a quasi-inessential saddle separatrix and a quasi-essential saddle separatrix (resp. two quasi-essential saddle separatrices).

\begin{figure}[t]
\begin{center}
\includegraphics[bb=0 0 102 102,scale=0.49]{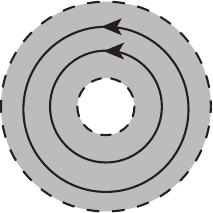}
\end{center}
\caption{Two-dimensional structure: a periodic annulus.}
\label{fig:pericdic_annulus}
\end{figure}

As a two-dimensional structure, we also introduce an open annulus filled with periodic orbits, called a {\it periodic annulus}, as shown in Figure~\ref{fig:pericdic_annulus}.

\subsection{Topological structure theorem for Hamiltonian flows on a flat tours}~\label{sec:3.2}
 A $C^r$-Hamiltonian vector field $X$ ($r\geq 1$) on a compact surface is said to be {\it structurally stable} if the topological structure of the streamlines does not change when small disturbances are added to the vector field. 
 The exact definition is given below.
\begin{definition}
Let $\mathcal{H}^r(M)$ be the set of $C^r$-Hamiltonian vector fields ($r\geq 1$) on a compact surface $M$. 
A $C^r$-Hamiltonian vector field $X$ on $M$ is said to be {\it structurally stable} (in $\mathcal{H}^r(M)$) if any Hamiltonian vector field $\widetilde{X}$ that is $C^1$-close to $X$ is topologically equivalent to $X$.
That is, there exists a homeomorphism $h \colon M \to M$ such that $h$ maps the orbits of $X$ to the orbit of $\widetilde{X}$ without changing the direction of the orbit. 
\end{definition}

The following fact about structurally stable Hamiltonian vector fields on compact surfaces is known.
\begin{theorem}{\cite[Theorem 2.3.8, p. 74]{ma2005geometric}}\label{th:str_stable}
Let $X$ be a $C^r$ Hamiltonian vector field on a compact surface. The flow generated by $X$ is denoted by $v_X$.
Then the following are equivalent:
\begin{enumerate}
\item $X$ is structurally stable.
\item Any singular orbits of $v_X$ are non-degenerate and the saddle connection diagram of $v_X$ is self-connected.
\end{enumerate}
Furthermore, the set of structurally stable $C^r$-Hamiltonian vector fields is an open and dense subset of the set of $C^r$-Hamiltonian vector fields.
\end{theorem}

In this paper, with an application to two-dimensional turbulence in mind, we assume that the flow domain does not contain any physical boundaries.
Then we have the following lemma.

\begin{table}[]
 \begin{center}
\begin{tabular}{|c|c|c|}\hline
$0$-dim. structures & saddles, centers & Figure~\ref{fig:Point_Flow_Components} \\ \hline
\multirow{2}{*}{$1$-dim. structures} & self-connected saddle separatrices, & Figure~\ref{fig:ess_saddle_connection}\\
& periodic boundary components & \\ \hline
$2$-dim. structures & open periodic annuli & Figure~\ref{fig:pericdic_annulus} \\ \hline
\end{tabular}
\end{center}
\caption{Orbit structures appearing in structurally stable Hamiltonian flows on a torus.}
\label{tbl:fundamental}
\end{table}

\begin{figure}[t]
\begin{center}
\includegraphics[scale=0.45]{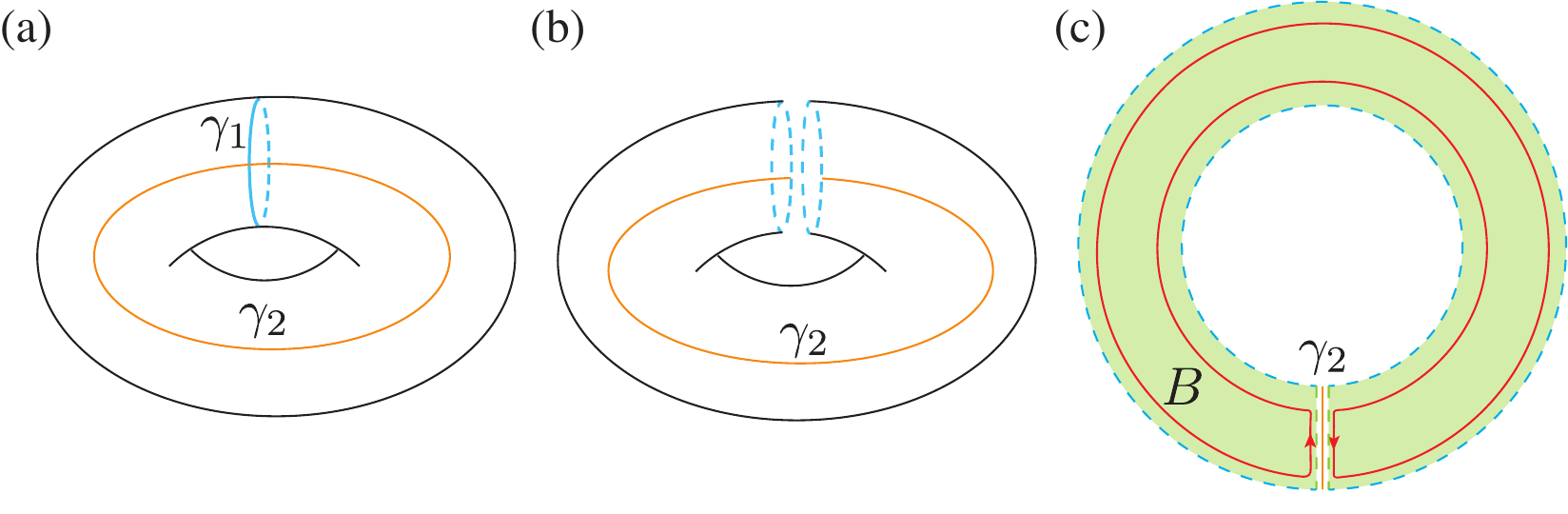}
\end{center}
\caption{(a) A fully essential saddle connection $C$ on a torus. (b) We remove the essential simple closed curve $\gamma_1$ from the torus. (c) The open annulus $A=\mathbb{T}^2-\gamma_1$ and the open disk $B=A \setminus C$ are filled with green.}
\label{fig:annulus_separatrix}
\end{figure}

\begin{lemma}\label{lem:non-fully}
Every self-connected saddle connection in a Hamiltonian flow on a torus without a physical boundary is not fully essential.
\end{lemma}

\begin{proof}
Let $C$ be a self-connected saddle connection of the Hamiltonian flow on the torus $\mathbb{T}^2$. 
Then it is a union of the saddle and its two connecting separatrices.

We prove it by contradiction. 
Assume that $C$ is fully essential as shown in Figure~\ref{fig:annulus_separatrix}(a). 
Then we choose a saddle separatrix whose closure is an essential simple closed curve $\gamma_1 \subset C$. 
Removing $\gamma_1$ from the torus, we have an open annulus $A:= \mathbb{T}^2 - \gamma_1$. 
Note that the intersection $A \cap C$ is the other separatrix, say $\gamma_2$. 
See Figure~\ref{fig:annulus_separatrix}(b).
Therefore, the set difference $B:= A \setminus C = \mathbb{T}^2 -C$ is an invariant open disk as shown in Figure~\ref{fig:annulus_separatrix}(c). 
Hence $C$ is the boundary $\partial B$, and so is a boundary component of a periodic annulus. 

On the other hand, let us consider a periodic orbit near $\partial B = \gamma_1 \sqcup \gamma_2$, which appears as a closed red orbit in Figure~\ref{fig:annulus_separatrix}(c). 
Then, the flow directions of the periodic orbit along both sides of $\gamma_2$ should be the opposite. However, the flow follows the curve $\gamma_2$ in one direction, which is a contradiction. 

Since any saddle connection of a structurally stable flow on a torus without a physical boundary is self-connected, the second assertion follows from the first assertion. 
\end{proof}
It follows from Theorem~\ref{th:str_stable} and Lemma~\ref{lem:non-fully} that orbit structures in structurally stable Hamiltonian flows on a torus without a physical boundary are topologically characterized as follows.

\begin{proposition}\label{lem:ssc}
The following holds for structurally stable Hamiltonian flows on a torus without a physical boundary:
\begin{enumerate}
\item Each self-connected saddle connection is one of the four one-dimensional structures in Figure~\ref{fig:ess_saddle_connection} up to isotopy on the torus.
\item The singular orbits of the Hamiltonian flow are either centers or saddles.
\end{enumerate}
\end{proposition}

\begin{proof}
By Theorem~\ref{th:str_stable} and \cite[Theorem 3]{cobo2010flows}, the second assertion holds.
Let $C$ be a self-connected saddle connection. 
By Lemma~\ref{lem:non-fully}, the self-connected saddle connection $C$ is not fully essential. 
If $C$ is inessential, then $C$ is one of the four one-dimensional structures in Figure~\ref{fig:ess_saddle_connection}(a) up to isotopy on the torus.
If $C$ is singly essential, then $C$ is one of them in Figure~\ref{fig:ess_saddle_connection}(b,c) up to isotopy on the torus.
\end{proof}
This proposition indicates that the streamlines of structurally stable instantaneous Hamiltonian are decomposed into a finite number of orbit structures given in Table~\ref{tbl:fundamental}.

\section{Tree representation of topological orbit structures}\label{sec:4}
\subsection{Reeb graph for structurally stable Hamiltonian flows}\label{sec:4.1}
In TFDA, we convert the topological structures of structurally stable Hamiltonian flows on a torus into a discrete plane tree, called a partially Cyclically Ordered rooted Tree (COT for short).
To this end, we first introduce the concept of a \textit{finite graph} as follows.
\begin{definition}
A compact Hausdorff space $G$ is said to be a finite graph if it is a finite one-dimensional cell complex.
In other words, there exists a finite number of points $v_1, \ldots, v_k \in G$ such that $G - \{v_1, \ldots, v_k\}$
 is homeomorphic to a finite disjoint union of open intervals.
\end{definition}
For a real-valued function $f \colon X \to \R$ on a topological space $X$ and $c \in \R$, the inverse image $f^{-1}(c)$ is called a {\it level set}.
Furthermore, for a function $f \colon X \to \R$, the space obtained by collapsing the connected components of each level set into a singleton is called the Reeb graph. 
The space is denoted by $X/f$.
More precisely, when there exists a connected component $C$ of $f^{-1}(c)$ for some $c \in \R$ such that $p, q \in C$, we define the equivalence relation $p \sim q$. 
Then the Reeb graph is represented by $X/f:= X/\sim$.
In what follows, we consider the Reeb graph of the Hamiltonian for a structurally stable Hamiltonian flow.

\begin{lemma}\label{lem:Reeb_finite}
The Reeb graph of a Hamiltonian on a compact surface $S$ for a structurally stable Hamiltonian flow is a finite graph.
\end{lemma}

\begin{proof}
Let $v$ be a structurally stable Hamiltonian flow generated by a Hamiltonian $H$ on a compact surface. 
By Proposition~\ref{lem:ssc}, the singular point of $v$ is either a saddle or a center.
Let $D(v)$ be the saddle connection diagram. 
The finiteness of singular points implies that $D(v)$ consists of finitely many orbits. 
Therefore, the complement of them in $S$ consists of a finite number of periodic annuli.
As a result, the nodes and edges of the Reeb graph correspond to these singular points and the periodic annuli, respectively.
\end{proof}
This construction method allows us to associate a Hamiltonian value with each node in the Reeb graph, which corresponds to a saddle point or a center.
It is also possible to assign the difference in Hamiltonian values of the nodes at both ends to each edge as an attribute. 
We also note that each node for a center has degree one since each center corresponds to a local maximum or minimum of the Hamiltonian. 
On the other hand, since a saddle divides the domain into three regions by its self-connected saddle connections, as shown in Figure~\ref{fig:ess_saddle_connection}, the node corresponding to the saddle has degree three.
This means that the Reeb graph of the Hamiltonian of a structurally Hamiltonian flow is a finite graph consisting of nodes of degree $1$ and $3$.

\subsection{Essential periodic orbits}\label{sec:4.2}
\begin{figure}[t]
\begin{center}
\includegraphics[scale=0.27]{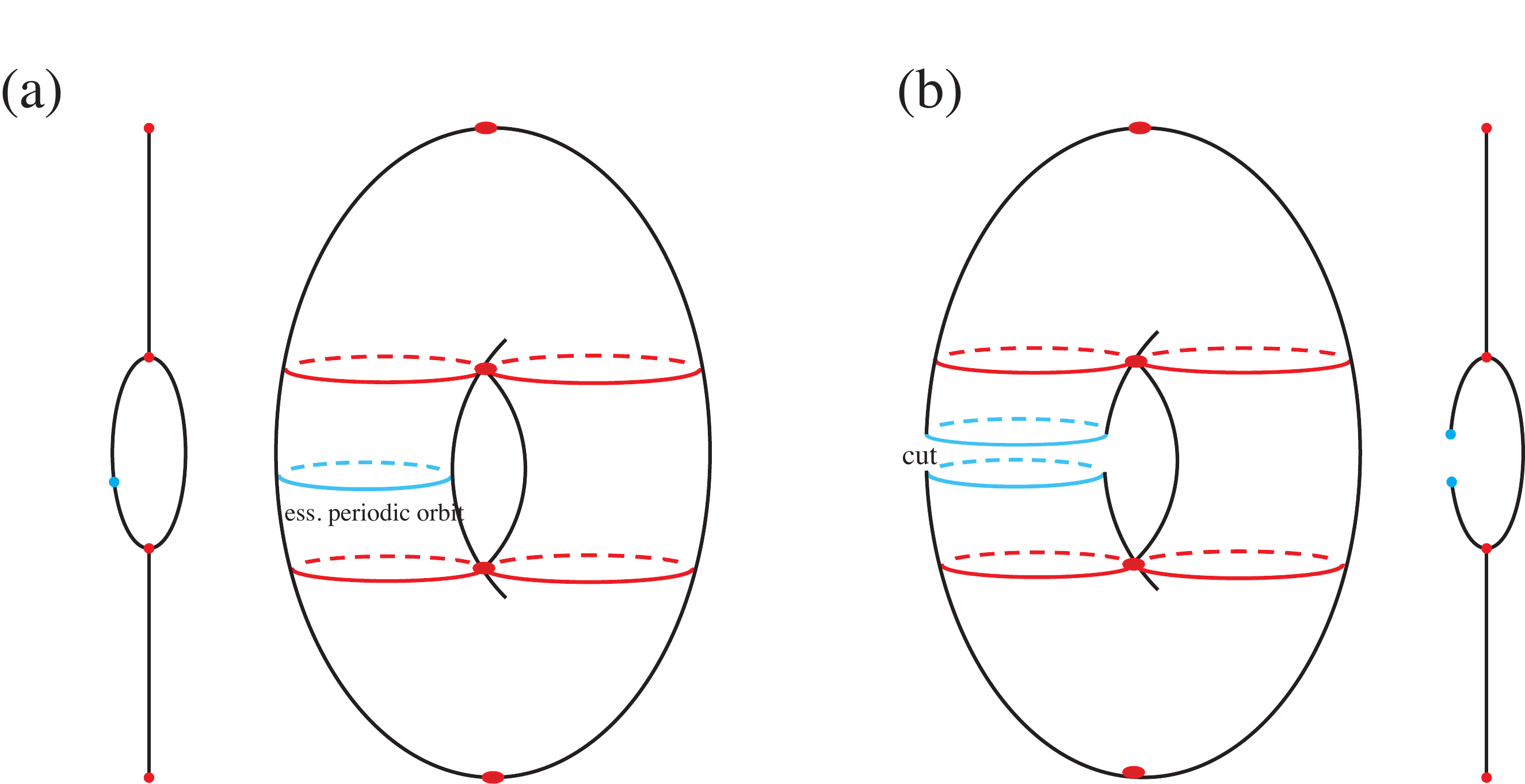}
\end{center}
\caption{(a) The Reeb graph for a Hamiltonian on the torus $\mathbb{T}^2$. 
The Reeb graph of the function has a loop structure since the Hamiltonian is split into two at the handle structure of the torus.  
(b) When we cut one of the handles along an essential periodic orbit, the Reeb graph of 
the Hamiltonian becomes a tree.}
\label{fig:Reeb_graph}
\end{figure}

Figure~\ref{fig:Reeb_graph}(a) shows the Reeb graph for a monotonically increasing Hamiltonian on a torus, with the minimum point at the bottom of the torus and the maximum point at the top.
As shown in this figure, the Reeb graph of a continuous function on a torus has a loop structure and is not a tree in general.
Since the output of TFDA is a tree, we need to address this situation. 
To this end, we introduce {\it an essential periodic orbit}, which is characterized as follows.
\begin{lemma}
For a periodic orbit $O$ of a Hamiltonian flow on a torus, the following conditions are equivalent:
\begin{enumerate}
\item The periodic orbit $O$ is essential.
\item The periodic orbit $O$ does not bound a disk.
\item The complement $\mathbb{T}^2 - O$ is an annulus.
\end{enumerate}
\end{lemma}

\begin{proof}
Any sufficiently small $r$-neighborhood of $O$ is homotopic to $O$.
If $O$ is essential, then $\mathbb{T}^2 - O$ is an annulus by definition. 
If $\mathbb{T}^2 - O$ is an annulus, then $O$ does not bound a disk.
If $O$ does not bound a disk, then $O$ is essential.
\end{proof}

For a torus without a physical boundary, we can show the existence of essential periodic orbits for structurally stable Hamiltonian flows on the torus as follows.

\begin{lemma}\label{lem:ess.per}
Every structurally stable Hamiltonian flow on a torus without a physical boundary has essential periodic orbits.
\end{lemma}

\begin{proof}
Let $v$ be a structurally stable Hamiltonian flow on a torus and $H$ be a Hamiltonian that generates the flow $v$. 
By Proposition~\ref{lem:ssc}, each singular orbit is either a saddle or a center, and each self-connected saddle connection is one of the four one-dimensional structures in Figure~\ref{fig:ess_saddle_connection} up to the isotopy on the torus.
Hence, every self-connected saddle connection is inessential or singly essential. 
From Poincar{\'e}--Hopf theorem, since the Euler characteristic of the torus is zero and the index of a saddle (resp. a center) is $-1$ (resp. $1$), the number of saddles equals that of the centers. 
From Lemma~\ref{lem:Reeb_finite}, the Reeb graph $G$ of the Hamiltonian $H$ is a connected finite graph. 
Moreover, the saddles correspond to trivalent nodes, and the centers correspond to univalent nodes.

We claim that the Reeb graph $G$ is not a tree. 
In fact, assume that $G$ is a tree. 
Since saddles correspond to trivalent (degree-three) nodes and the centers to univalent (degree-one) nodes by the nonexistence of loops in the tree $G$, the number of centers is at least two more than the number of saddles, which contradicts the equality of the numbers of centers and saddles.

Since $G$ is not a tree, the graph $G$ contains a loop $K \subset G$. 
Then any periodic annuli corresponding to the edges in the loop $K$ are essential. 
Therefore, so are any periodic orbits contained in the periodic annuli. 
\end{proof}
\begin{figure}[t]
\begin{center}
\includegraphics[width=.85\linewidth]{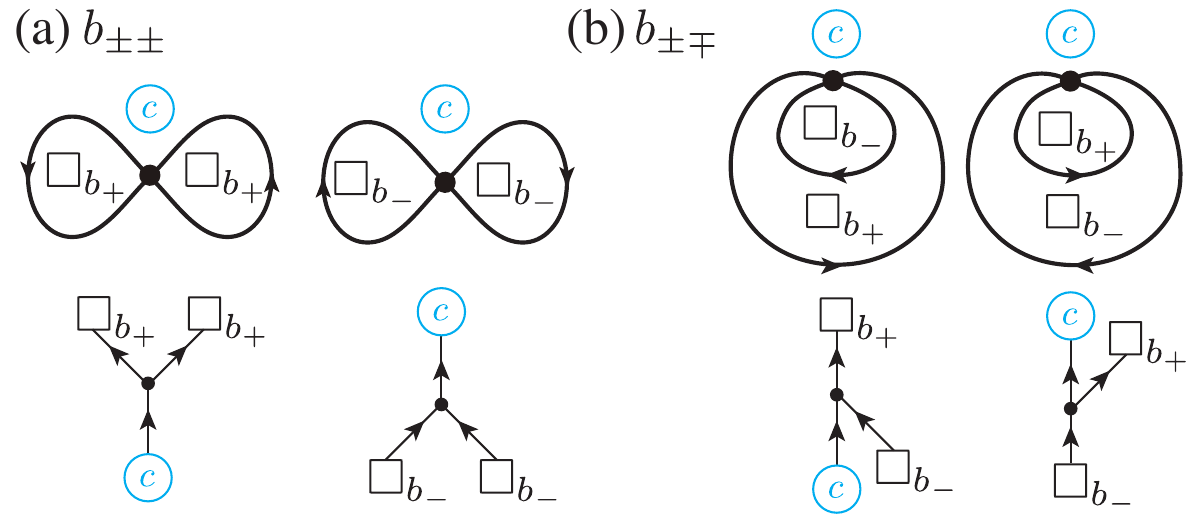}
\end{center}
\caption{One-dimensional orbit structures with a saddle and two quasi-inessential saddle separatrices. 
(a) $b_{\pm\pm}$ structures: figure-eight orbit structures.  Their COT symbols are given by $b_{\pm\pm} \{ \B_{b_\pm}, \B_{b_\pm} \}$. 
(b) $b_{\pm\mp}$ structures: orbit structures where one quasi-inessential saddle separatrix encloses another quasi-inessential saddle separatrix. Their COT symbols are given by $b_{\pm\mp}(\B_{b_\pm}, \B_{b_\mp})$.}
\label{fig:S2_structures}
\end{figure}

\subsection{Root and current edge}\label{sec:4.3}
The Reeb graph $G$ of the Hamiltonian for a structurally stable Hamiltonian flow becomes a finite graph with exactly one loop, in general, as shown in Figure~\ref{fig:Reeb_graph}(a). 
On the other hand, Lemma~\ref{lem:ess.per} guarantees the existence of essential periodic orbits in structurally stable Hamiltonian flows on a torus without a physical boundary. 
Hence, by cutting the torus along one of the essential periodic orbits, we can always make the flow domain an annular one whose boundary is the essential periodic orbit.
This operation corresponds to cutting one edge of the loop from the Reeb graph $G$ into two disjoint edges, producing a tree $G'$, as shown in Figure~\ref{fig:Reeb_graph}(b). 
Then the cut periodic orbit corresponds to the two nodes of the tree. 
We thus choose the root node as the lower one, which is a local maximum of the Hamiltonian in the periodic annulus.
Notice that the cutting operation increases by adding one node and one edge each.
Unless otherwise confused, when the root node has degree one, the union of the root and its connecting edges is also referred to as \textit{root}. 

To explain the conversion algorithm to COT in Section~\ref{sec:4.5}, we define the notion of a {\it current} for the edges of a rooted tree.
Let $G =(V, E)$ be a tree with a root node $v_0 \in V$. 
Fix a non-root node $v \in V - \{ v_0 \}$. 
An edge $e$ of $G$ is a current of $v$, if $v$ is contained in the boundary of $e$ and the minimal path from the root to $v$ contains $e$. 
For the rooted tree $G'$ of the Reeb graph $G$ as above, a maximal invariant open periodic annulus that corresponds to the current of a non-root node is called a current domain.

\subsection{Local orbit structures and COT symbols}\label{sec:4.4}
As discussed in Section~\ref{sec:4.3}, the Reeb graph of the Hamiltonian becomes a tree by cutting an edge at the point corresponding to an essential periodic orbit. 
A partially Cyclically Ordered rooted Tree (COT) is then obtained by assigning a label, which specifies the topological orbit structure of a saddle or a center, to each node of this tree.
To this end, we define the following strings,  called COT symbols, for each orbit structure.

\bigskip
\noindent
\textbf{($\sigma_{\pm}$ structures)} The simplest zero-dimensional local orbit structure is a center in Figure~\ref{fig:Point_Flow_Components}(a).
We assign the COT symbol $\sigma_+$ (resp. $\sigma_-$) to a center associated with counterclockwise (resp. clockwise) periodic orbits in its neighborhood, where the subscript $\pm$ expresses their flow direction.
They are the label of nodes with degree $1$ in the Reeb graph since each center is a local maximum or minimum of the Hamiltonian.

\bigskip
\noindent
\textbf{($b_{\pm\pm}$, $b_{\pm\mp}$ structures)} They are one-dimensional orbit structures consisting of a saddle and two quasi-inessential self-connected separatrices, as shown in Figure~\ref{fig:ess_saddle_connection}(a).

For the figure-eight structure in Figure \ref{fig:S2_structures}(a), we assign the COT symbol $b_{++}\{\B_{b_+}, \B_{b_+}\}$ (resp. $b_{--}\{\B_{b_-}, \B_{b_-}\}$) when the rotational direction of the two self-connected saddle separatrices is counterclockwise (resp. clockwise). 
Each saddle separatrix encloses an inner structure symbolized by $\B_{b_\pm}$, chosen from one of the following local orbit structures.
\begin{equation}
\B_{b_+} \in \{ b_{++}, b_{+-}, \sigma_+\}, \qquad \B_{b_-}  \in \{ b_{--}, b_{-+}, \sigma_-\}. 
\label{b-structures}
\end{equation}
Note that the rotational direction of the inner structures is determined by that of the separatrices. 
In addition, the arrangement order of the two inner structures is cyclic, which is expressed by the curly braces in the COT symbols.

The COT symbol $b_{\pm\mp}(\B_{b_\pm}, \B_{b_\mp})$ is provided for the orbit structures in Figure \ref{fig:S2_structures} (b), where a quasi-inessential saddle separatrix encloses another quasi-inessential saddle separatrix.  
When the outer self-connected saddle separatrix is going counterclockwise (resp. clockwise) direction, the COT symbol is given by $b_{+-}(\B_{b_+}, \B_{b_-})$ (resp. $b_{-+}(\B_{b_-}, \B_{b_+})$), in which $\B_{b_\pm}$ denotes the inner structure chosen from (\ref{b-structures}) and the round braces denote that the arrangement of the local orbit structures is uniquely determined.

The node of degree $3$ in the Reeb graph of the Hamiltonian for these two orbit structures is also shown in Figure~\ref{fig:S2_structures}(a,b).

\begin{figure}[tp]
\begin{center}
\includegraphics[scale=0.45]{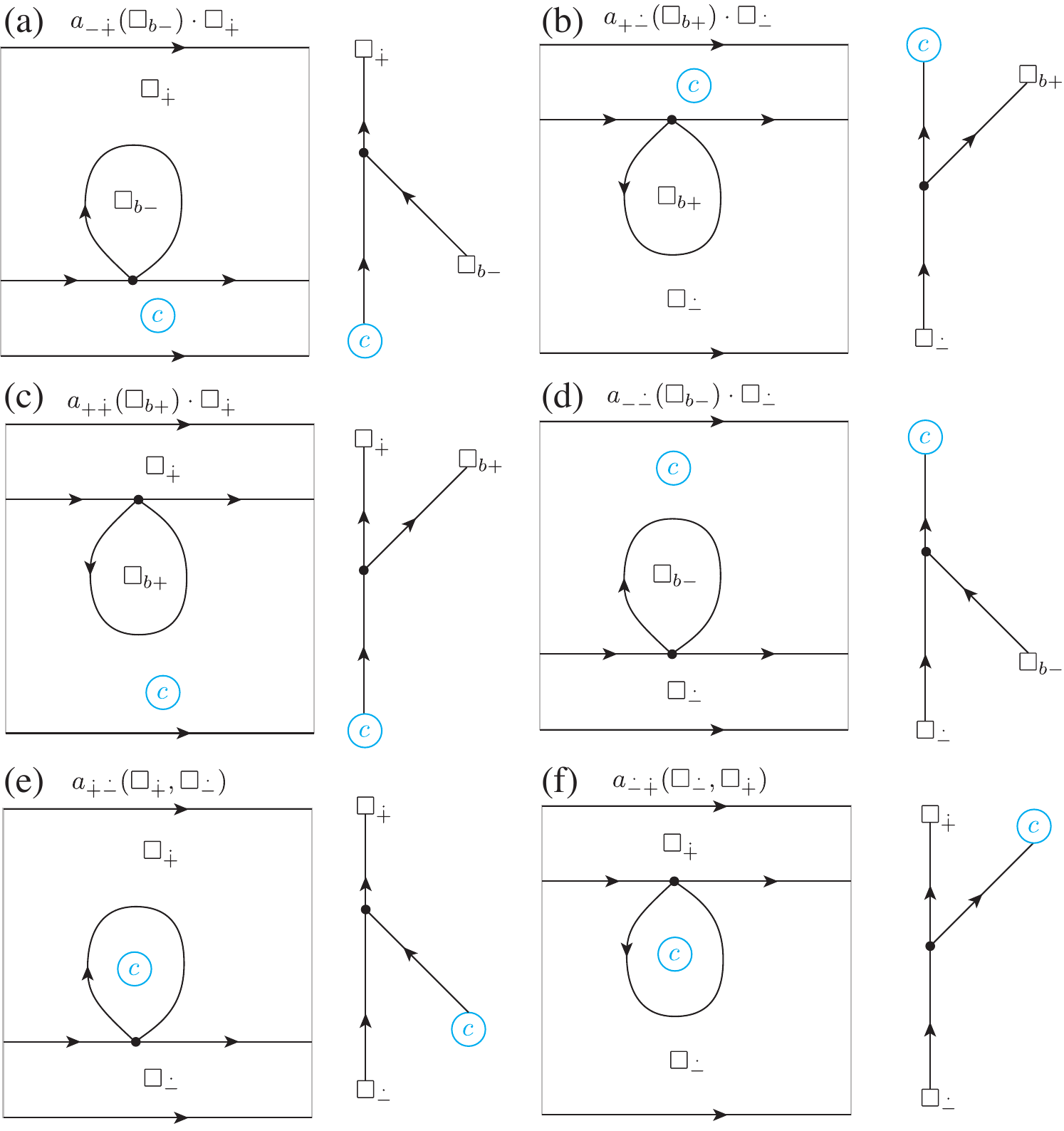}
\end{center} 
\caption{One-dimensional structures consisting of one quasi-essential self-connected saddle separatrix and one quasi-inessential saddle connection. 
For each structure, the local structure of a node of degree $3$ in the Reeb graph of the Hamiltonian is shown in each panel.
The current domain is denoted by $\copyright$.}
\label{cot01}
\end{figure}

\bigskip
\noindent
\textbf{($a_{\pm \dot{+}}$, $a_{\pm\dot{-}}$, $a_{\dot{+}\dot{-}}$, $a_{\dot{-}\dot{+}}$ structures)} We consider one-dimensional structures consisting of a quasi-inessential saddle separatrix and a quasi-essential saddle separatrix in Figure~\ref{fig:ess_saddle_connection}(b).
These saddle separatrices divide the flow domain into three.
In this case, six different structures must be identified as topologically distinct structures depending on the choices of the three current regions and the two flow directions. 
For each structure, the corresponding local structure of a node of degree $3$ in the Reeb graph of the Hamiltonian is shown in each panel of Figure~\ref{cot01}.
See Figure~\ref{cot01}(a)--(f).

\begin{itemize}
\item Suppose that the intersection of the boundary of the essential current domain $\copyright$ and the saddle separatrices does not contain the quasi-inessential saddle separatrix.
When the Hamiltonian value of the saddle separatrix is higher than those of the current domain, as shown in Figure~\ref{cot01}(a), we assign the COT symbol $a_{-\dot{+}}(\B_{b_-}) \cdot \B_{\dot{+}}$ to the orbit structure.
The symbol $\B_{b_-}$, chosen from (\ref{b-structures}), denotes a local orbit structure enclosed by the quasi-inessential saddle separatrix. 
The symbol $\B_{\dot+}$ represents a local orbit structure in the other essential domain, which is chosen from 
\begin{equation}
\Box_{\dot{+}} \in \{ a_{\pm \dot{+}}, \alpha_{+ \dot{-}}, \lambda_{\dot{+}} \}. \label{pm-structures+}
\end{equation}

On the other hand, the Hamiltonian value of the saddle separatrices is lower than those of the current domain, we have the  orbit structure in Figure~\ref{cot01}(b).
Its COT symbol is given by the COT symbol $a_{+\dot{-}}(\B_{b_+}) \cdot \B_{\dot{-}}$, in which $\B_{\dot-}$ is chosen from 
 \begin{equation}
\Box_{\dot{-}} \in \{  a_{\pm\dot{-}}, a_{-\dot{+}}, \lambda_{\dot{-}}\}. \label{pm-structures-}
\end{equation}
In these COT symbols, the first subscript $+$ (resp. $-$) of $a_{\pm\dot\mp}$ denotes that the flow direction of the quasi-inessential saddle separatrix is counterclockwise (resp. clockwise) and the second subscript $\dot+$ (resp. $\dot-$) indicates that the Hamiltonian values in the other essential domain are higher (resp. lower) than those in the current essential domain.

\item Suppose that the intersection of the boundary of the essential current domain $\copyright$ and the saddle separatrices contains the quasi-inessential saddle separatrix. 
When the Hamiltonian value of the saddle separatrices is higher (resp. lower) than those of the current domain, we have the orbit structure in Figure~\ref{cot01}(c) (resp. Figure~\ref{cot01}(d)) and its COT symbol is given by $a_{+\dot{+}}(\B_{b_+}) \cdot \B_{\dot{+}}$ (resp. $a_{-\dot{-}}(\B_{b_-}) \cdot \B_{\dot{-}}$), in which $\Box_{b_\pm}$ and $\Box_{\dot\pm}$ are chosen from (\ref{b-structures}) and (\ref{pm-structures+})--(\ref{pm-structures-}) respectively.
\item Suppose that the current domain $\copyright$ is inessential. 
When the Hamiltonian value of the saddle separatrices is higher (resp. lower) than those in the current domain, where the orbit structure is shown as Figure~\ref{cot01}(e) (resp. Figure~\ref{cot01}(f)), its COT symbol is given by $a_{\dot{+}\dot{-}}(\B_{\dot{+}},\B_{\dot{-}})$ (resp. $a_{\dot{-}\dot{+}}(\B_{\dot{-}},\B_{\dot{+}})$). 

The first subscript $\dot+$ (resp. $\dot-$) of $a_{\dot\pm\dot\mp}$ indicates that the Hamiltonian values of the essential domain that share the quasi-inessential saddle separatrix with the current domain $\copyright$ are higher (resp. lower) than the Hamiltonian value of the saddle.
The second subscript $\dot-$ (resp. $\dot+$) denotes that the Hamiltonian values of the other essential domain are lower (resp. higher) than the Hamiltonian value of the saddle.
\end{itemize}
\begin{figure}[t]
\begin{center}
\includegraphics[scale=0.45]{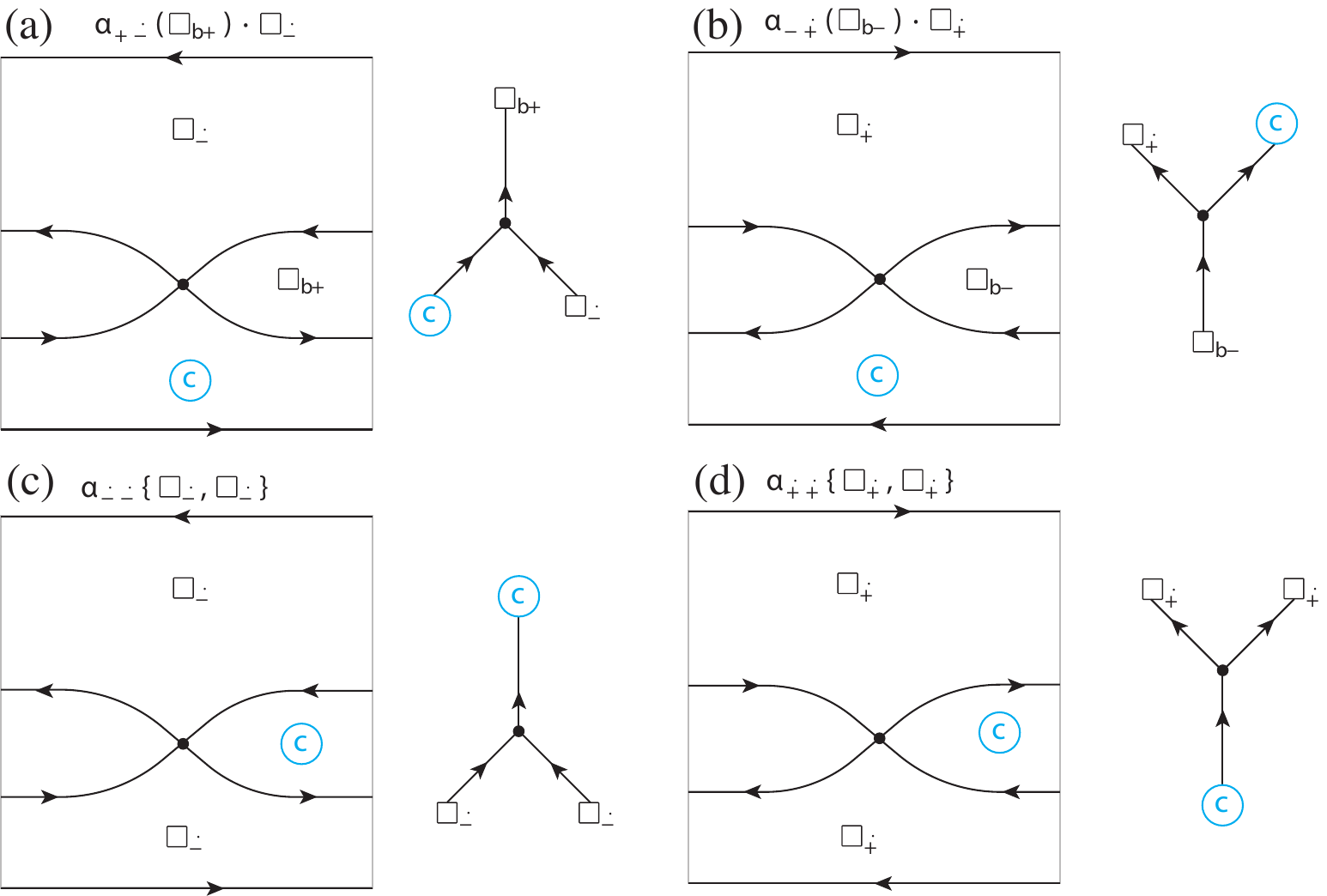}
\end{center} 
\caption{One-dimensional structures consisting of two quasi-essential self-connected saddle separatrices. 
For each structure, the local structure of a node of degree $3$ in the Reeb graph of the Hamiltonian is shown in each panel.} The current domain is denoted by $\copyright$.
\label{cot01+}
\end{figure}

\noindent
\textbf{($\alpha_{\dot{+}\dot{+}}$, $\alpha_{\dot{-}\dot{-}}$, $\alpha_{- \dot{+}}$, $\alpha_{+ \dot{-}}$ structures)} We consider orbit structures with two quasi-essential saddle separatrices in Figure~\ref{fig:ess_saddle_connection}(c).
The quasi-essential separatrices divide the flow domain into three parts, but due to the rotational symmetry of the flow structure, the current domain can be set in one of the two essential domains or in the essential domain. 
Consequently, we have four topologically distinct orbit structures, as shown in Figure~\ref{cot01+}(a)--(d).
For each structure, the corresponding local structure of a node of degree $3$ in the Reeb graph of the Hamiltonian is shown in each panel of Figure~\ref{cot01+}.

\begin{itemize}
\item Suppose that the current domain $\copyright$ is essential. 
When the Hamiltonian values of the current domain are lower (resp. higher) than that of the saddle connection, we have the orbit structure of Figure~\ref{cot01+}(a) (resp. Figure~\ref{cot01+}(b)). 
The COT symbol is given by $\alpha_{+ \dot{-}}( \B_{b_+}) \cdot \B_{\dot{-}}$ (resp. $\alpha_{- \dot{+}}( \B_{b_-}) \cdot \B_{\dot{+}}$), in which $\Box_{b_\pm}$ and $\Box_{\dot\pm}$ are chosen from (\ref{b-structures}) and (\ref{pm-structures+})--(\ref{pm-structures-}) respectively.

\item Suppose that the current domain $\copyright$ is inessential. 
When the Hamiltonian values of the current domain are higher (resp. lower) than that of the saddle connection, we have the orbit structure of Figure~\ref{cot01+}(c) (resp. Figure~\ref{cot01+}(d)). The COT symbol is given by $\alpha_{\dot{-}\dot{-}}\{ \B_{\dot{-}}, \B_{\dot{-}} \}$ (resp. $\alpha_{\dot{+}\dot{+}}\{ \B_{\dot{+}}, \B_{\dot{+}} \}$), in which the arrangement of the orbit structures in $\Box_{\dot\pm}$, chosen from (\ref{pm-structures+}) and (\ref{pm-structures-}) is cyclic due to rotational symmetry.
\end{itemize}

\noindent
\textbf{($\beta_{\dot\pm}$ structures)}  This one-dimensional orbit structure represents a boundary component of the annular domain that appears when the torus is cut by an essential periodic orbit.
We assign the COT symbol $\beta_{\dot+}$ (resp. $\beta_{\dot-}$) to the boundary components of the domain with the lower (resp. higher) Hamiltonian values as viewed from the essential periodic orbit.
They are also the label of nodes of degree $1$ in the Reeb graph of the Hamiltonian.
By definition, the COT necessarily has a pair of $\beta_{\dot{+}}$ and $\beta_{\dot{-}}$. 
Furthermore, by identifying the essential periodic orbits $\beta_{\dot{+}}$ and $\beta_{\dot{-}}$, the Reeb graph of the Hamiltonian flow on the torus with one cycle is recovered. 

\subsection{Conversion algorithm to COT} \label{sec:4.5}
Let $ \widehat{H} \colon \R^2 \to \R $ be a $C^1$-function on $\R^2$, which is $2\pi$-doubly periodic. 
That is, for any $x,y \in \R^2$ and $m,n \in \Z$, $\widehat{H}(x+2\pi m, y+2\pi n) = \widehat{H}(x,y)$ is satisfied. 
Then a Hamiltonian $H \colon \mathbb{T}^2 \to \R$ is induced on the torus $\mathbb{T}^2 := \R^2 / (2\pi\Z)^2$. 
Suppose that the local coordinate system is induced from $\R^2$ before taking the quotient set.
Then it determines the origin of the torus for a given data set in the context of data analysis. 
We set the origin at the bottom left point of the coordinate system, and define the $x$-axis as the horizontal direction and the $y$-axis as the vertical direction.
 
Suppose now that a Hamiltonian flow on a torus without a physical boundary is structurally stable.
Then, we can choose an essential periodic orbit whose existence is guaranteed by Lemma~\ref{lem:ess.per}.
When we cut the torus along this orbit, the Hamiltonian flow on the torus is transformed into a Hamiltonian flow on a closed annulus whose boundaries are the cut periodic orbits.
For structurally stable Hamiltonian flow on the annulus, the existing conversion algorithm in~\cite{uda2019persistent_en} is available.
The algorithm is explicitly described as follows.

\begin{description}
\item{\textbf{(Step 1)}} We consider all sub-level sets that contain local minimums of the Hamiltonian. 
We then find the boundaries of the sub-level sets containing an essential simple closed curve, which has a minimal Hamiltonian value among such essential simple closed curves. 
 \item{\textbf{(Step 2)}} Since there exist essential periodic orbits near the essential simple closed curve chosen by (Step 1) on both sides of the level set, we uniquely choose an essential periodic orbit according to the procedure (P) described in~\ref{sec:appA}.
 \item{\textbf{(Step 3)}} Cut the torus along this essential periodic orbit $O$ and transform it into an annular region. 
\item{\textbf{(Step 4)}} We consider the domain in which the Hamiltonian values are lower than the Hamiltonian value of $O$. 
 It is set as the initial current domain whose boundary is $O$, to which we assign the COT symbol $\beta_{\dot+}$.
Starting from the current domain, we identify a local orbital structure corresponding to the other boundary of this current domain.
The current is then moved to the domain represented by the box symbols of the COT symbol for this local orbit structure, and a local orbit structure with a saddle is identified for each box symbol.
\item{{\textbf{(Step 5)}}} By moving the current domain, we similarly identify all orbit structures with a saddle corresponding to the box symbols of the COT symbol in a recursive manner until the identified orbit structure becomes either a center with a COT symbol $\sigma_\pm$ or an essential periodic orbit with a COT symbol $\beta_{\dot-}$.
\end{description}
Each time we identify an orbit structure with a saddle in Step 4 and Step 5, we connect its corresponding node associated with the Hamiltonian values.
As a result, we obtain a Reeb tree of the Hamiltonian on the annulus whose nodes are labeled with the COT symbol. 
This is the COT as required.
Note that from this construction method, the COT label of the terminal nodes is $\sigma_\pm$ or $\beta_{\dot\pm}$.
In addition, by embedding the COT symbol of the identified local structure in the box symbols, a string called a COT representation is constructed.
We will explain how the conversion procedure works with examples in Section~\ref{sec:5.1}.
As long as the selection rule of an essential periodic orbit is fixed as in Step 2, the COT and its COT representation are uniquely obtained, but it depends on the rule in general.

Finally, we show that the first local orbit structure with a saddle identified in Step 4 is always represented by the COT symbol $\alpha_{-\dot{+}}$.

\begin{lemma}\label{lem:a-+}
Every simple essential closed curve chosen by the algorithm contains a quasi-essential saddle separatrix whose corresponding COT symbol is $\alpha_{-\dot{+}}$. 
\end{lemma}

\begin{proof}
Let $H$ be a Hamiltonian that generates a structurally stable Hamiltonian flow on a torus.
In Step 1, the algorithm chooses a closed sub-level set $\Gamma$ that contains the local minimums of the Hamiltonian and saddles not in its interior but on the boundary.
Then the boundary component $\gamma$ in $\Gamma$ is the local minimum in the essential components connected to the level set near $\gamma$. 
By construction of $\Gamma$, the set of sub-levels $\Gamma$ contains a local minimum but no local maximum. 
This implies that the quasi-essential separatrix $\mu$ in $\gamma$ appears in Figure~\ref{cot01}(a,d,e) and Figure~\ref{cot01+}(b,d). 

Assume that the quasi-essential separatrix $\mu$ in $\gamma$ appears as in Figure~\ref{cot01}(a,d,e).
We then consider the subdomain consisting of the essential level sets there.
The Hamiltonian values decrease monotonically in the subdomain from top to bottom. 
Hence, any values of the Hamiltonian in the subdomain are not minimal. 
This implies that the quasi-essential saddle separatrix $\mu$ is not minimal in the subdomain, which contradicts the local minimality of $\gamma$ containing the quasi-essential saddle separatrix $\mu$. 
Thus, the quasi-essential separatrix $\mu$ is contained as in Figure~\ref{cot01+}(b,d).

Finally, in Figure~\ref{cot01+}(b,d), since $\gamma$ is the local minimum, Step 3 makes a current domain essential in the original surface $S$. 
Therefore, the one-dimensional structure $\alpha_{-\dot{+}}$ of Figure~\ref{cot01+}(b) only satisfies the locally minimal condition of $H$ and the position of the current simultaneously.
\end{proof}

\section{Topological vortex structures in turbulent flows}\label{sec:5}

\subsection{Application of the conversion algorithm to 2D turbulence patterns.}\label{sec:5.1}
We demonstrate how the conversion algorithm is applied to free-decaying turbulence and enstrophy-cascade turbulence in a doubly periodic domain. 
Recall that the stream function is the Hamiltonian of structurally stable Hamiltonian flows.

\begin{figure}[ht]
\begin{center}
\includegraphics[scale=0.45]{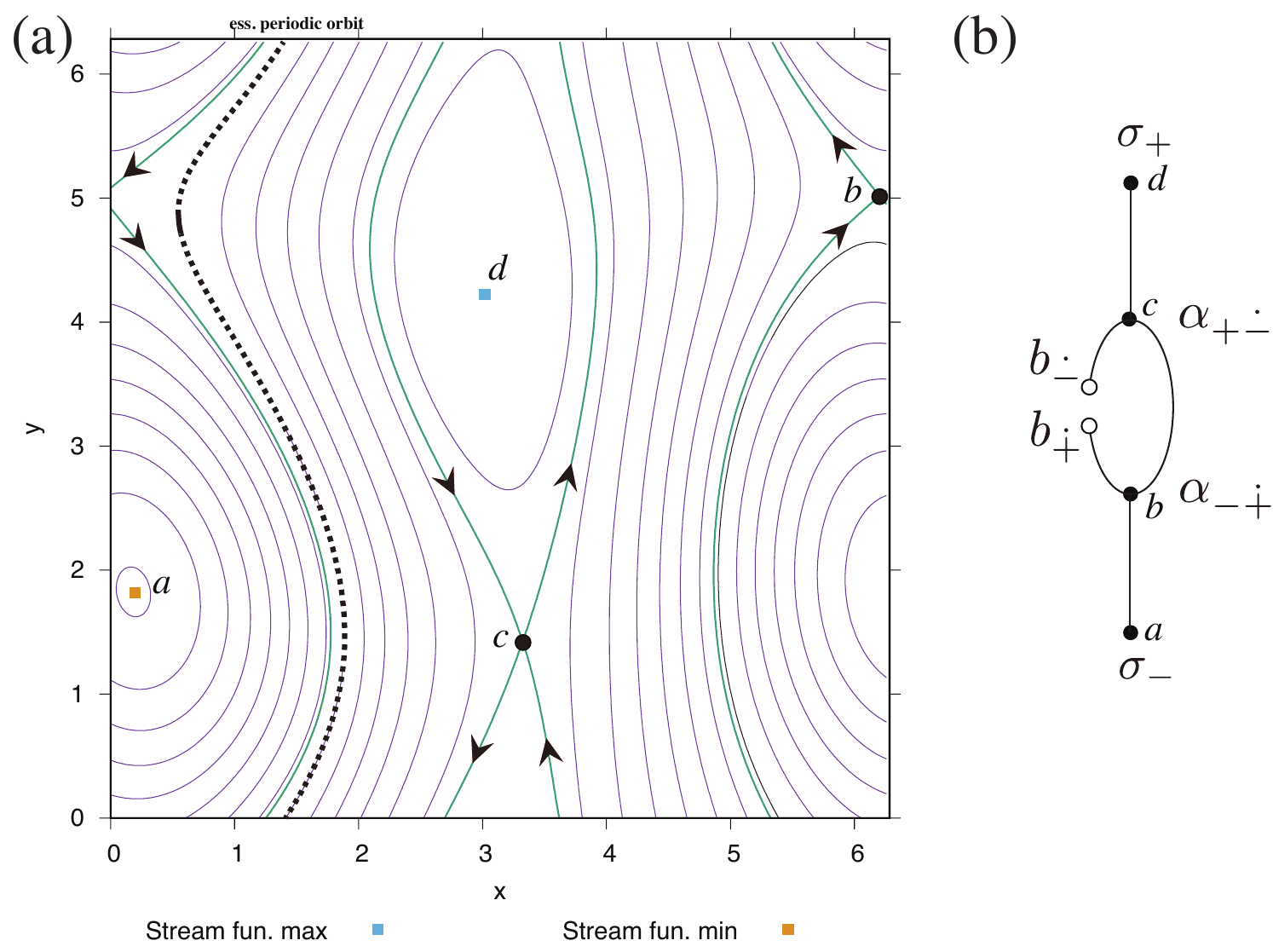}
\end{center} 
\caption{(a) A snapshot of level curves of the stream function of free decaying turbulence in a doubly periodic domain. 
The dashed level curve is the essential periodic orbit extracted by the procedure (P). 
(b) The COT of the Hamiltonian flow with the stream function
being its Hamiltonian. 
Its COT representation is (\ref{cot-rep-free}).}
\label{decay_2d_turb}
\end{figure}

Figure~\ref{decay_2d_turb}(a) shows a snapshot of the stream function for free-decaying turbulence.
Following Step 1 of the algorithm, we find the point where the Hamiltonian becomes minimum, which is attained in the center $a$ shown as the orange square in the figure. 
We then consider the domain that contains the minimum point between the two quasi-essential saddle separatrices of the saddle $b$.  
In Step 2, we choose an essential periodic orbit according to the procedure (P) in \ref{sec:appA}, which is shown as a dashed curve to the right of the saddle $b$.
In Step 3, we cut the torus along the dashed curve and obtain two boundaries of the annular domain. 

We proceed to Step 4. 
Since the Hamiltonian value in the region to the left of this essential periodic orbit decreases, we set the initial current in this region and provide the COT symbol $\beta_{\dot+}$ to the cut periodic orbit.
Starting from the current, we construct the COT and its COT representation.
Since the hierarchy of Hamiltonian values between the current domain and the quasi-essential separatrices of the saddle $b$ corresponds to the situation in Figure~\ref{cot01+}(b), we assign the COT symbol $\alpha_{-\dot+}(\B^1_{b_-})\cdot \B^1_{\dot+s}$ to this structure, as guaranteed by Lemma~\ref{lem:a-+}.
Moving to Step 5, we have $\B^1_{b_-}=\sigma_-$ since the local orbit structure for the box symbol $\B^1_{b_-}$ of $\alpha_{-\dot+}$ is the elliptic center $a$, and the process ends for this structure.
We then let the current move to the domain, which is an essential domain to the left of the separatrices of the saddle $b$, to identify the orbit structure in $\B^1_{\dot+s}$.
For this current, we find the saddle $c$ with two quasi-essential saddle separatrices as the local orbit structure to be identified. 
From the hierarchy of Hamiltonian values between the current domain and the saddle separatrices of the saddle $c$, we find that the local orbit structure is topologically equivalent to that of Figure~\ref{cot01+}(a). 
Hence, we have $\B^1_{\dot+s}=\alpha_{+\dot-}(\B^2_{b_+}) \cdot \B^2_{\dot-s}$ in which the inner structure is the elliptic center $d$ and thus $\B^2_{b_+}=\sigma_+$. 
Finally, the region to the left of the saddle $c$ has the essential periodic orbit that is the end of the structure. 
We finally assign $\B^2_{\dot-s}=\beta_{\dot-}$ to the essential boundary periodic orbit. 
Consequently, the COT is shown in Figure~\ref{decay_2d_turb}(b) whose COT representation is given by
\begin{equation}
\beta_{\dot+} \cdot \alpha_{-\dot+}(\sigma_-) \cdot \alpha_{+\dot-}(\sigma_+) \cdot \beta_{\dot-}. \label{cot-rep-free}
\end{equation}

\begin{figure}[ht]
\begin{center}
\includegraphics[scale=0.7]{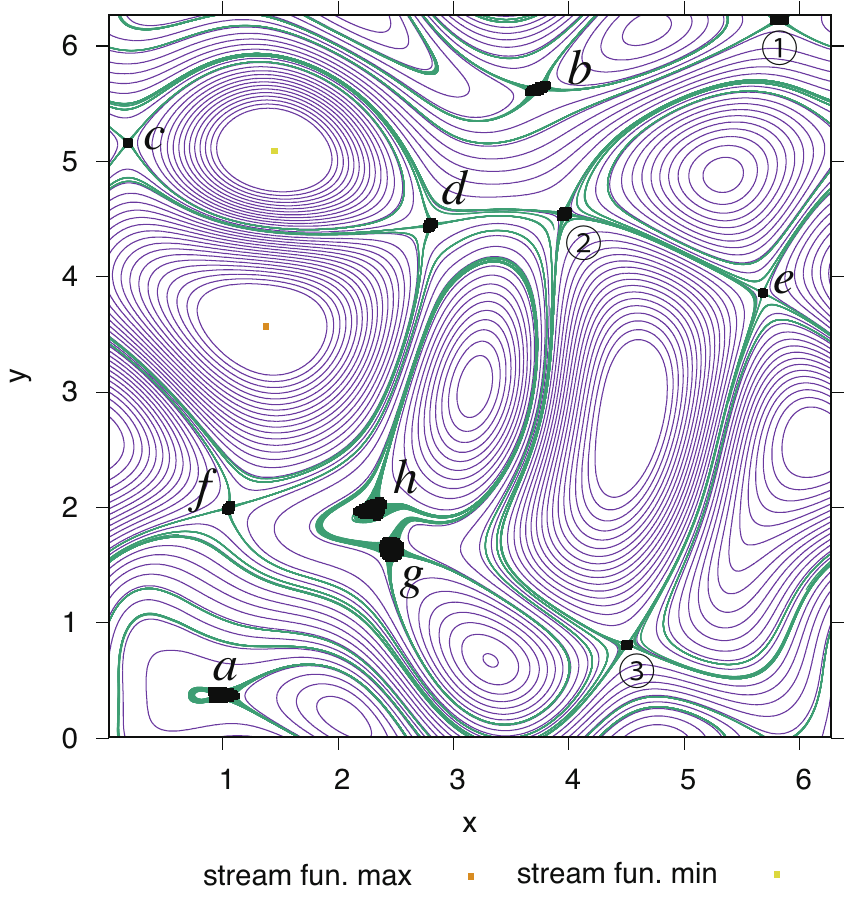}
\end{center} 
\caption{A snapshot of the level curves of the stream function of two-dimensional enstrophy cascade turbulence.}
\label{enstrophy_cascade_2d_turb}
\end{figure}

\begin{figure}[ht]
\begin{center}
\includegraphics[scale=0.5]{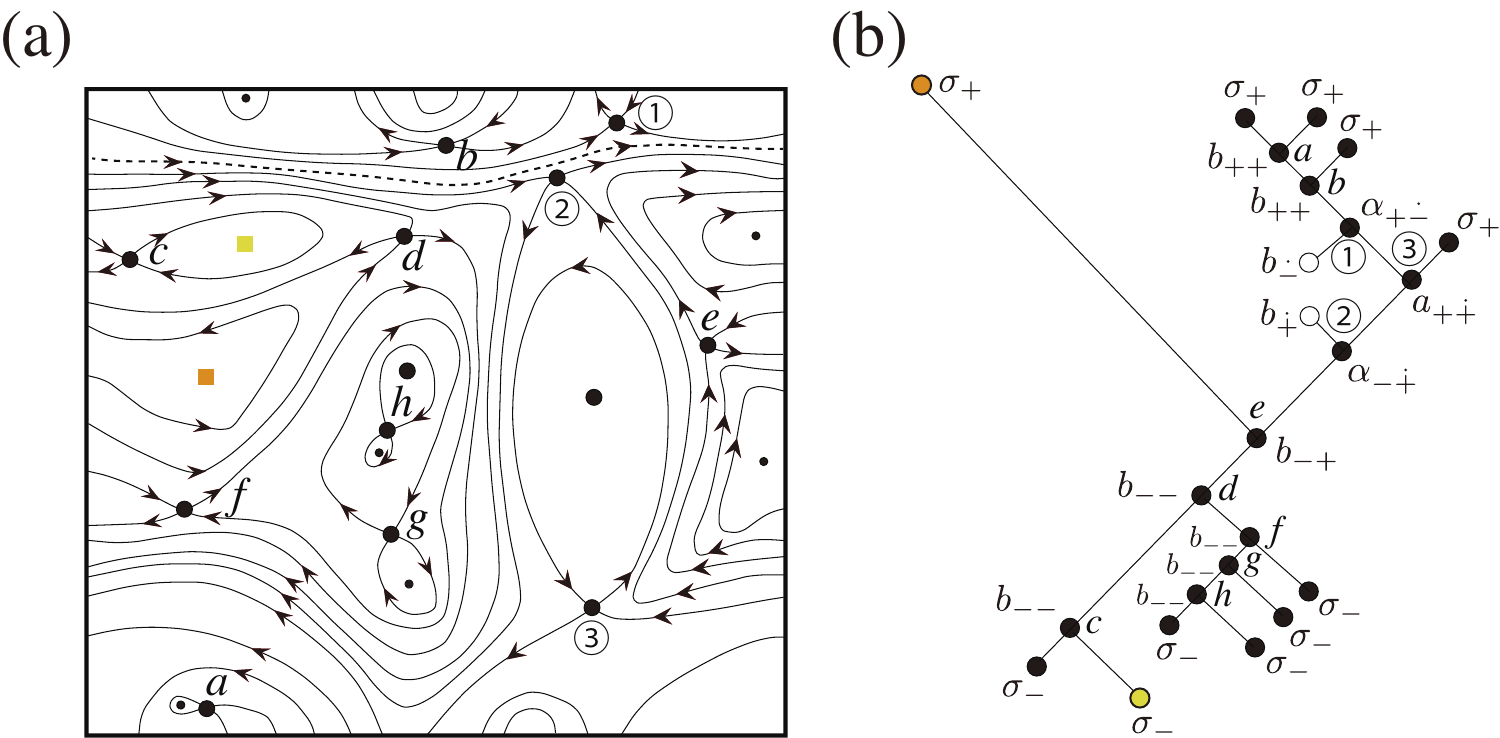}
\end{center} 
\caption{(a) The saddle connection diagram extracted schematically from the level curves of the stream function in Figure~\ref{enstrophy_cascade_2d_turb}.
(b) The COT for the saddle connection diagram. Its COT representation is given by (\ref{cot-rep-enstrophy}).}
\label{enstrophy_cascade_2d_turb_COT}
\end{figure}

Another example is an instantaneous stream function for enstrophy-cascade turbulence in Figure~\ref{enstrophy_cascade_2d_turb}. 
The yellow square (resp. the orange square) on the left of the figure represents the minimum point (resp. the maximum point) of the stream function.
This flow has three saddles (\ding{172}--\ding{174}) with quasi-essential saddle separatrices and seven saddles ($a$--$g$) with quasi-inessential self-connected saddle separatrices. 
Since the saddle separatrices are too densely packed in a small area to be distinguished, we redraw the saddle connection diagram schematically in Figure~\ref{enstrophy_cascade_2d_turb_COT}(a) to make it easy to apply the conversion algorithm. 

First, searching sub-level sets from the minimum point to higher Hamiltonian values, we find the saddle \ding{173} with quasi-essential saddle separatrices, which is to be identified in Step 1.
Following the procedure (P), we choose an essential periodic orbit in the upper region, which is drawn as a dashed curve, and cut the region along the orbit, which corresponds to Step 2 and Step 3.
We assign the COT symbol $\beta_{\dot+}$ to the essential periodic orbit and set the initial current domain in the region with lower Hamiltonian values below. 
The next local orbit structure we find for this current domain is the saddle \ding{173} with two quasi-essential saddle separatrices whose configuration is topologically equivalent to that in Figure~\ref{cot01+}(h). 
Hence, the COT symbol for this structure is given by $\alpha_{-\dot+}(\B^1_{b_-}) \cdot \B^1_{b\dot+}$.
The orbit structure in $\B^1_{b_-}$ is the saddle $e$ whose topological structure is represented by the COT symbol $\B_{b_-}^1 = b_{-+}(\B^2_{b_-}, \B^2_{b_+})$.
Since the inner structure $\B^2_{b_+}$ of $b_{-+}$ is the elliptic center having the maximum Hamiltonian value, we have $\B^2_{b_+}=\sigma_+$ and terminate the process.
For the other inner structure in $\B^2_{b_-}$ of $b_{-+}$, we find the saddle $d$ represented by the COT symbol $b_{--}\{\B^3_{b_-}, \B^4_{b_-}\}$. 
Similarly, checking the inner structure $\B^3_{b_-}$, we find a nested structure of the saddles $f$, $g$ and $h$, each of which is represented by the COT symbol $b_{--}$.
Hence, we have $\B^3_{b_-}=b_{--}\{b_{--}\{b_{--}\{\sigma_-,\sigma_-\}, \sigma_-\}, \sigma_-\}$.
On the other hand, for the inner structure $\B^4_{b_-}$, we have the saddle $c$ with the COT symbol $b_{--}$, which yields $\B^4_{b_-}=b_{--}\{\sigma_-,\sigma_-\}$.

Next, we move the current to the domain of $\B^1_{b\dot+}$ below the saddle \ding{173} to identify the other structure in $\B^1_{b_-}$ of $\alpha_{-\dot+}$. 
Searching the topological structure for higher values of the Hamiltonian, we find the saddle \ding{174} with a quasi-essential saddle separatrix and a counterclockwise quasi-inessential saddle separatrix whose configuration is topologically equivalent to that in Figure~\ref{cot01}(c). 
Hence, the COT symbol for this local orbit structure becomes $\B^1_{b\dot+}=a_{+\dot+}(\B^5_{b_+})\cdot \B^5_{\dot+}$. 
The inner orbit structure in $\B^5_{b_+}$ corresponds to the region containing a center $\sigma_+$. 
Hence, we set $\B^5_{b_+}=\sigma_+$, and the search ends. 
We let the current move to the other region of $\B^5_{\dot+}$, and search topological orbit structures toward higher values of the Hamiltonian. 
We find the saddle \ding{172} with quasi-essential saddle separatrices whose configuration is topologically equivalent to Figure~\ref{cot01+}(a) represented by the COT symbol $\B^5_{\dot+}=\alpha_{+\dot-}(\B^6_{b_+})\cdot \B^6_{\dot-}$.
For the structure $\B^6_{\dot-}$, we find the essential periodic orbit $\B^6_{\dot-}=\beta_{\dot-}$ and the search in this direction ends.
Finally, searching for the topological orbit structure towards higher values of the Hamiltonian for $\B^6_{b_+}$, we have a nested structure of the saddles $a$ and $b$ with the COT symbol $b_{++}$. 
Hence, we have $\B^6_{b_+}=b_{++}\{b_{++}\{\sigma_+,\sigma_+\},\sigma_+\}$.
Since all topological orbit structures have been detected, we have the COT for the Hamiltonian flows as in Figure~\ref{enstrophy_cascade_2d_turb_COT}(b).
The COT representation is finally given by
\begin{multline}
\beta_{\dot+} \cdot\alpha_{-\dot+}(b_{-+}(b_{--}\{b_{--}\{\sigma_-, \sigma_-\}, b_{--}\{b_{--}\{b_{--}\{\sigma_-,\sigma_-\}, \sigma_-\}, \sigma_-\}\},\sigma_+))  \\
\cdot a_{+\dot+}(\sigma_+) \cdot \alpha_{+\dot-}(b_{++}\{b_{++}\{\sigma_+, \sigma_+\},\sigma_+\})\cdot \beta_{\dot-}. \label{cot-rep-enstrophy}
\end{multline}

\subsection{Statistical study of topological vortex structures in turbulent patterns}\label{sec:5.2}
We apply our TFDA method to statistically steady 2-D numerical turbulent flows in a doubly periodic domain. 
The COT representations for these flows become much more complex.
Hence, we consider statistics of certain quantities extracted from TFDA.
It is known that the 2-D Navier--Stokes turbulence at high Reynolds numbers can take two different states: the one is the energy inverse cascade (IC) state and the other is the enstrophy cascade (EC) state; see, e.g., \cite{Tabeling2002-by}.
This is due to the presence of two conservative quantities of the Euler equations in two dimensions, namely the kinetic energy $\mathcal{E}$ and the enstrophy $\mathcal{L}$, which are defined by
\[
\mathcal{E}(t) = \frac{1}{2}\int_{\mathbb{T}^2} \vert {\bm u}({\bm x},t)\vert^2 d{\bm x}, \qquad 
\mathcal{L}(t)= \frac{1}{2}\int_{\mathbb{T}^2} \vert \omega({\bm x},t)\vert^2 d{\bm x}
\]
with the velocity field ${\bm u}({\bm x},t)$ and the vorticity $\omega({\bm x},t)$.
To let the flow reach a statistically steady state, we need to add a stirring force. 
The standard force used in laboratory or numerical experiments of two-dimensional turbulent flows is confined to a certain wavenumber $k_f$. 
It is known that the IC state is realized for the wavenumber range $k \ll k_f$ and the EC state is realized for $k \gg k_f$, provided that the size of the system is sufficiently large and the viscous dissipation is sufficiently small.
However, in practice, satisfying both requirements is demanding. 
Therefore, we often set the forcing scale $2\pi/k_f$ to be comparable to either the system size or the characteristic length scale of the viscous dissipation. 
The former setting yields only the EC state that spans nearly all wavenumbers and the latter yields only the IC states \cite{boffetta_two-dimensional_2012}.
By adopting this convention, we here numerically obtain the IC and EC states in separate simulations.
More specifically, we use the numerical methods in \cite{ecnum,icnum} for the IC and EC respectively using the $8$th order Laplacian for the viscosity to have a cleaner scaling property, the so-called hyper-viscosity.

\begin{figure}
    \centering
    \includegraphics[width=0.8\linewidth]{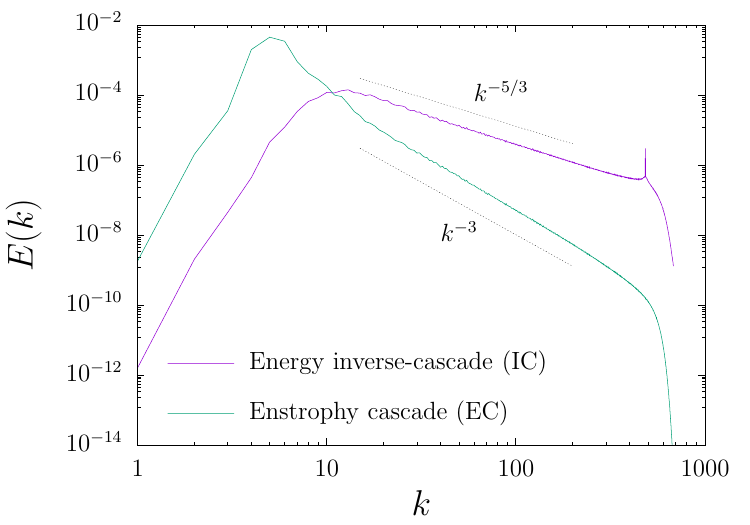}
    \caption{Time-averaged energy spectra of 2-D turbulence for the energy inverse-cascade (IC) state and for the enstrophy cascade (EC) state. The forcing terms for the IC state and for the EC state are added in $484.5 < |{\bm k}| < 485.5$ and in $4.0 < |{\bm k}| < 6.0$, respectively. 
    The forcing term added to the Navier--Stokes equations is proportional to the velocity in those wavenumber ranges and its prefactor is set to ensure the constant energy injection for the IC case and the constant enstrophy injection for the EC case. 
    For both states, a large-scale energy drag is added to the Navier-Stokes equations, which is proportional to $\nabla^{-2} {\bm u}$ to obtain a statistically steady state. For both cases, the number of grid points is $2048^2$ and the wavenumber bin-size $\Delta k$ for the energy spectrum is $1$. The periodic domain is square with a side length $2\pi$.
}
    \label{f:spc}
\end{figure}

In the conventional turbulence study, the time average of the energy spectrum $E(k)$ is used to distinguish between the IC and EC states, where $k$ is the wavenumber. 
The energy spectrum is defined by
\[
E(k) = \sum_{k \le \vert{\bm k}\vert < k + \Delta k} \frac{1}{2}\left\vert \widehat{\bm u}({\bm k}, t)\right\vert^2 \frac{1}{\Delta k},
\]
where $\widehat{{\bm u}}({\bm k}, t)$ denotes the Fourier coefficient of the wave vector ${\bm k}$ at time $t$ for the velocity field.
Figure~\ref{f:spc} shows the energy spectra of our simulated IC and EC states. 
They are consistent with the phenomenological predictions, $E(k) \propto \epsilon^{2/3} k^{-5/3}$ for the IC state and $E(k) \propto \eta^{2/3} k^{-3}$ for the EC state in their respective inertial ranges. 
Here, $\epsilon$ and $\eta$ are the $k$-independent energy flux and enstrophy flux, respectively.
Note that the more precise form of the EC spectrum is $E(k) \propto \eta^{2/3} k^{-3} \ln^{-1/3}(k / k_f)$, see, e.g., \cite{boffetta_two-dimensional_2012}. 
Note that we ignore the logarithmic correction of the EC state.

Using TFDA, we propose a method for vortex identification in these turbulent patterns from a topological perspective.
Although several other methods of extracting vortex structures are known \cite{rivera_direct_2014}, our identification method unveils several new statistical aspects on the nature of turbulent fluctuations in the IC and EC states.
Here, we describe details of our topological data analysis.
We prepare $2500$ snapshots for both evolutions in the statistically steady states, which are equally spaced with the time step of $0.2$ for time $50$.
For each instantaneous Hamiltonian, we first coarse-grain the Hamiltonian values on the $2048\times 2048 $ simulation grids to those on a $256 \times 256$ grid stored at each coarse-grained point (pixel) to reduce the size of the data.
We then normalize the coarse-grained Hamiltonian values.
Specifically, the Hamiltonian values at each time step are divided by the difference between the maximum and minimum values at that time step.
It has been shown in \cite{uda2019persistent_en} that even with such size reduction and normalization, the topological orbit structures of the Hamiltonian flow are substantially extracted by TFDA.

We need to perform TFDA for a large number of snapshots, and therefore manual topological data analysis conducted in Section~\ref{sec:5.1} is no longer practical.
To this end, we use $\mathbb{T}^2$-\texttt{psiclone}, a Python library that implements the TFDA algorithm for doubly periodic Hamiltonian flows developed in~\cite{sakamoto_tjsiam_2025submitted}.
This software not only performs TFDA quickly on large amounts of data but also provides a useful additional feature: a low-pass filter for the scale of orbit structures.
Recall that every edge is associated with the difference between the Hamiltonian values at the two end nodes as its weight, as explained in Section~\ref{sec:4.1}.
Hence, for a given threshold $\epsilon_0>0$, we can remove certain edges with weights less than or equal to $\epsilon_0$. This filtering is used later.

\begin{figure}[htbp]
    \centering
    \begin{subfigure}[b]{0.3\textwidth}
        \includegraphics[width=\textwidth]{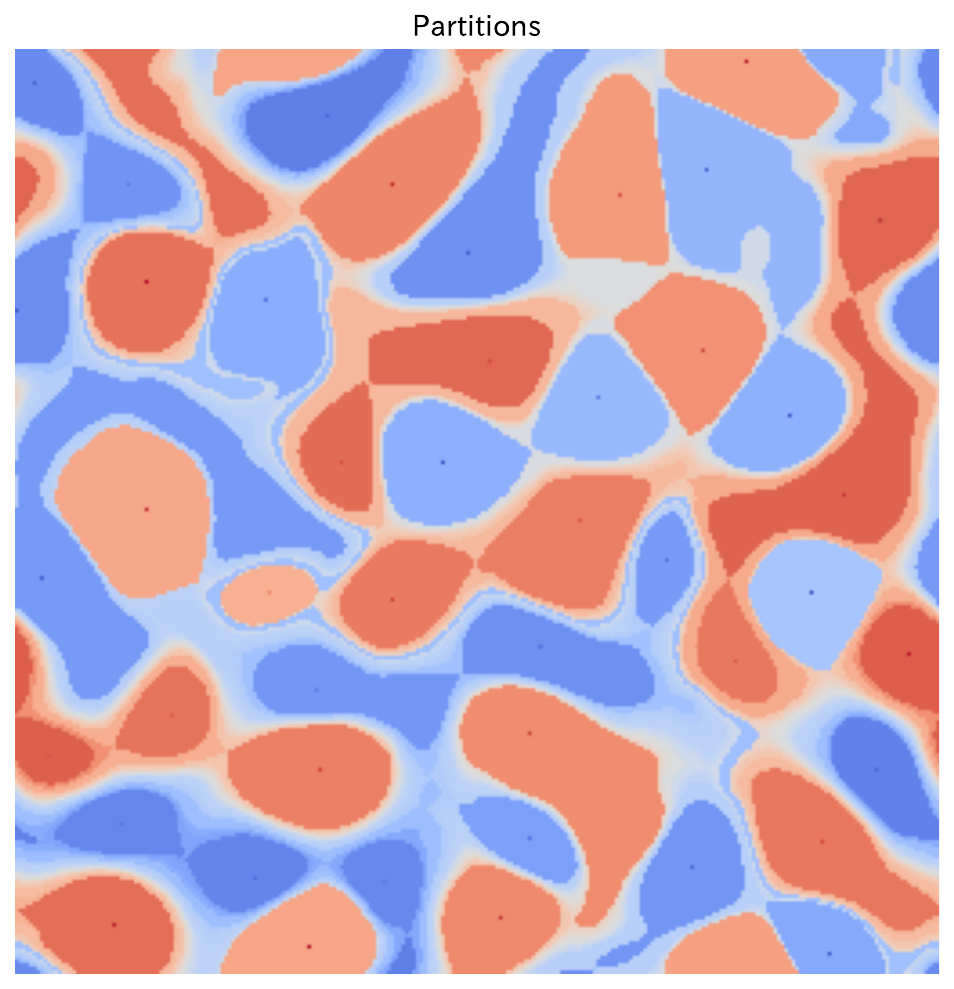}
        \caption{}
        \label{fig:terminal_vortex_all}
    \end{subfigure}
    \hfill
    \begin{subfigure}[b]{0.3\textwidth}
        \includegraphics[width=\textwidth]{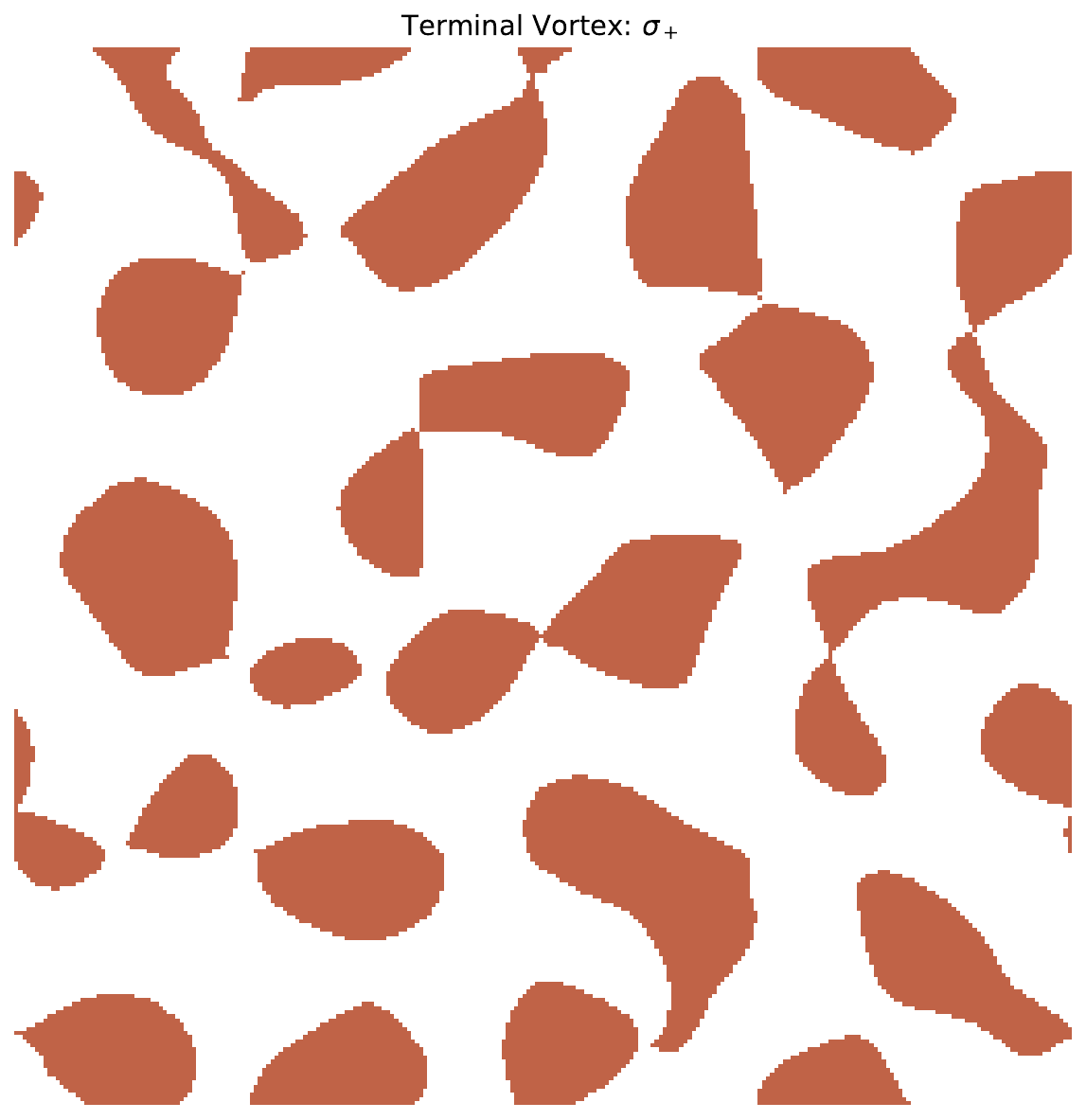}
        \caption{}
        \label{fig:terminal_vortex_plus}
    \end{subfigure}
    \hfill
    \begin{subfigure}[b]{0.3\textwidth}
        \includegraphics[width=\textwidth]{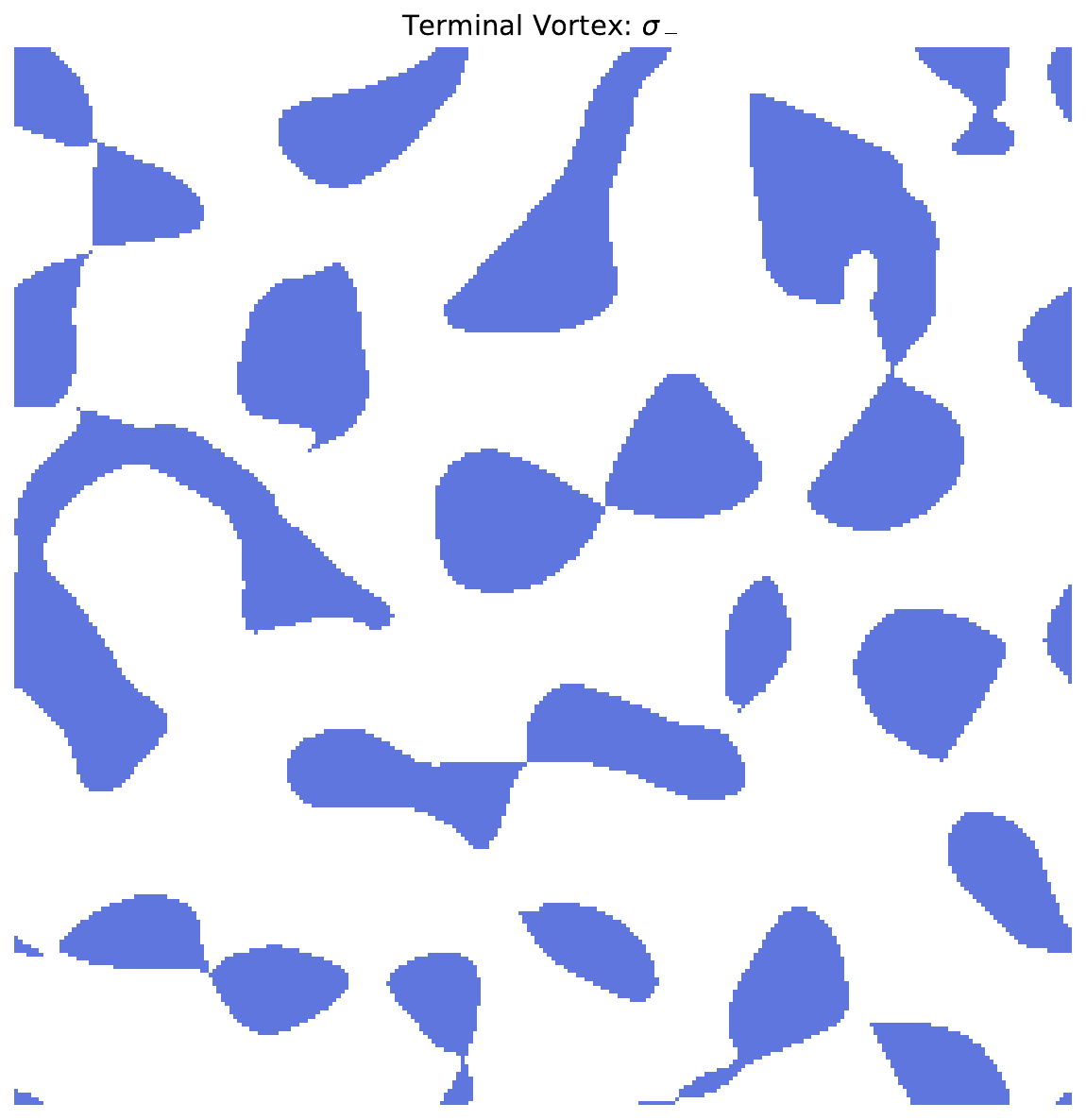}
        \caption{}
        \label{fig:terminal_vortex_minus}
    \end{subfigure}
 \caption{Visualization of terminal vortices for a turbulent flow data (a snapshot of the EC state). (a) The topological partitioning of the flow domain based on COT representation computed by TFDA. (b) Counterclockwise terminal vortices corresponding to $\sigma_+$. (c) Clockwise terminal vortices for $ \sigma_-$, respectively.}
    \label{fig:terminal_vortex}
\end{figure}

From the Reeb graph (COT) of the Hamiltonian, we introduce ``terminal vortices'' to identify vortex-like structures in 2D turbulent flows, which are the regions corresponding to terminal edges toward leaf nodes of the labels $\sigma_\pm$.
In this paper, using filtering with $\epsilon_0=0.1$, we remove all terminal edges whose weights do not exceed $\epsilon_0$, since they are considered topologically small structures.
Terminal vortices represent isolated swirling flows, as illustrated in Figure~\ref{fig:terminal_vortex} for a snapshot of the EC state. 
In general, we note that terminal vortices do not coincide with high-vorticity regions.
However, as observed in this figure, these terminal vortices are consistent with our intuitive picture of vortices.

\begin{figure}[htbp]
    \centering

    \begin{subfigure}[b]{0.48\textwidth}
        \includegraphics[width=\textwidth]{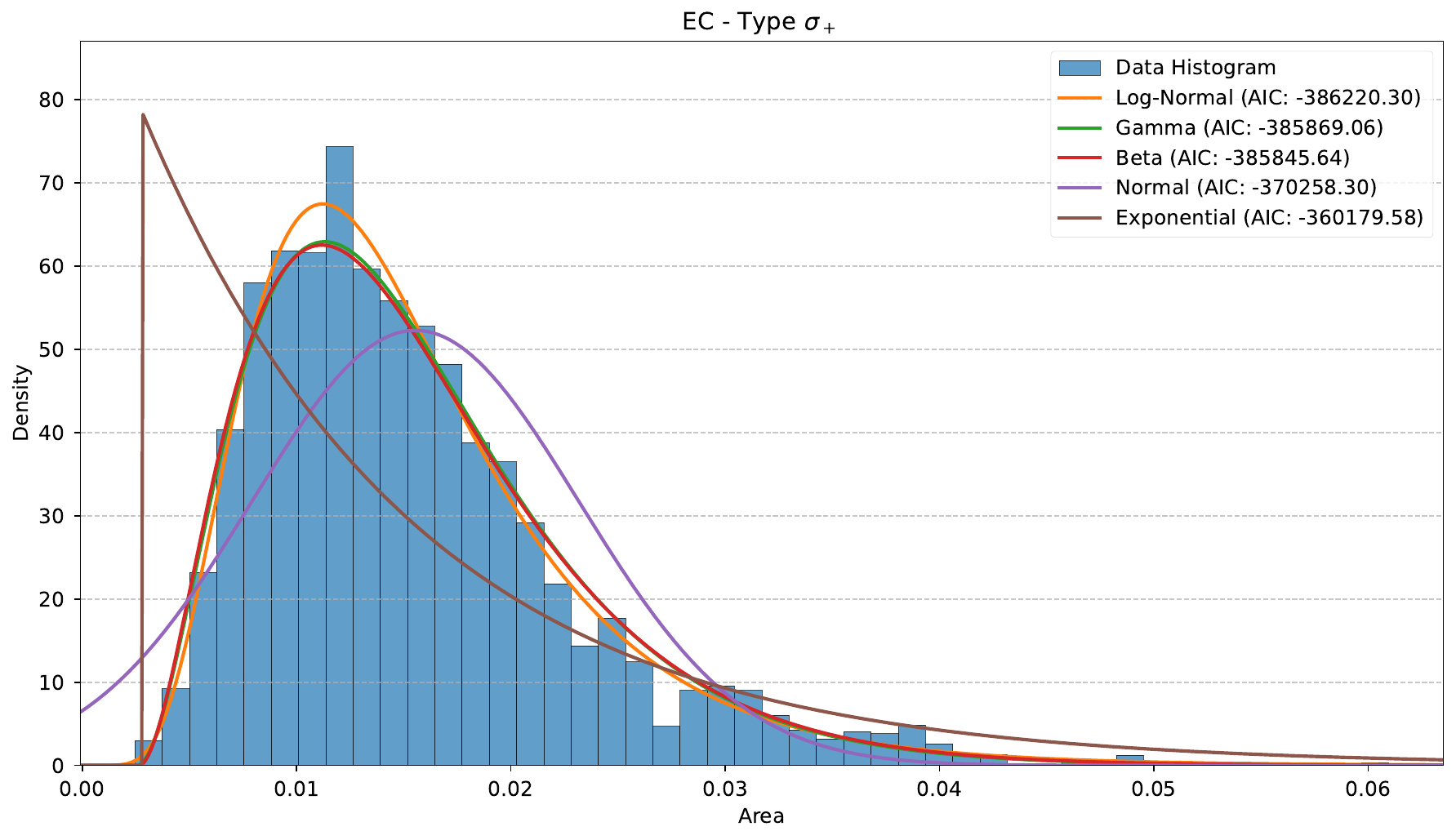}
        \caption{EC - Type $\sigma_+$}
        \label{fig:ec_sigma_plus_fits}
    \end{subfigure}
    \hfill
    \begin{subfigure}[b]{0.48\textwidth}
        \includegraphics[width=\textwidth]{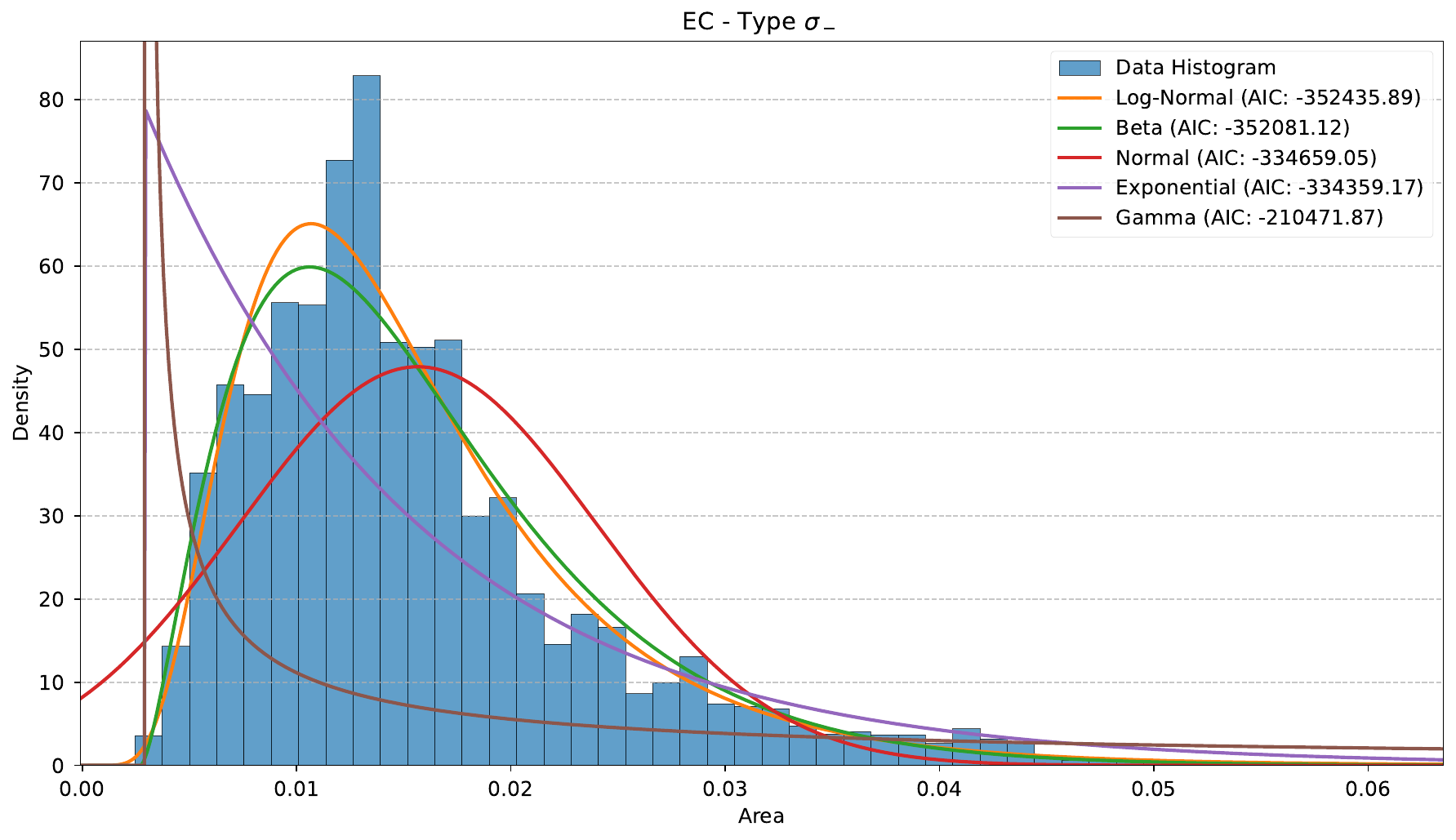}
        \caption{EC - Type $\sigma_-$}
        \label{fig:ec_sigma_minus_fits}
    \end{subfigure}

    \vspace{0.5cm}

    \begin{subfigure}[b]{0.48\textwidth}
        \includegraphics[width=\textwidth]{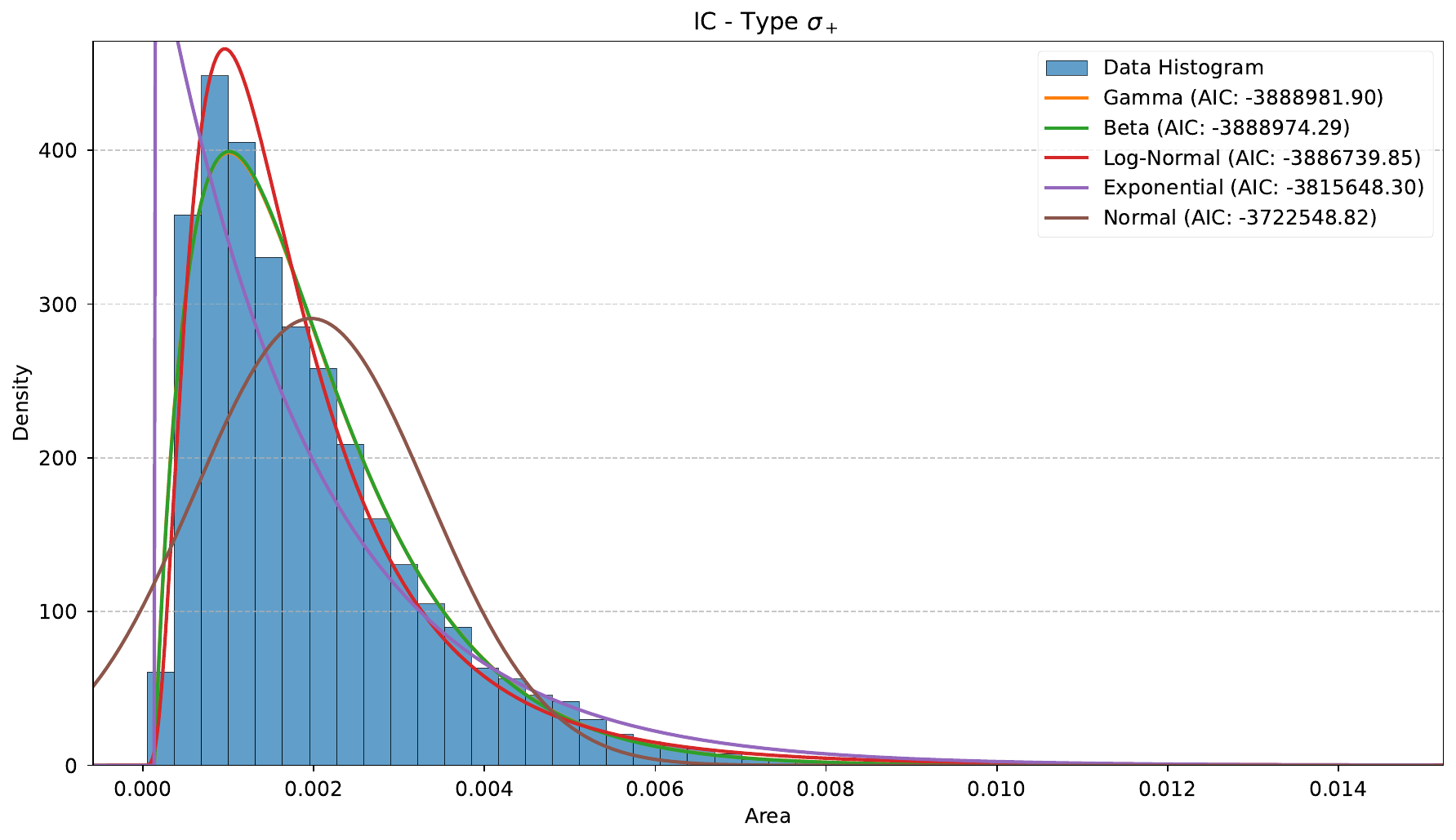}
        \caption{IC - Type $\sigma_+$}
        \label{fig:ic_sigma_plus_fits}
    \end{subfigure}
    \hfill
    \begin{subfigure}[b]{0.48\textwidth}
        \includegraphics[width=\textwidth]{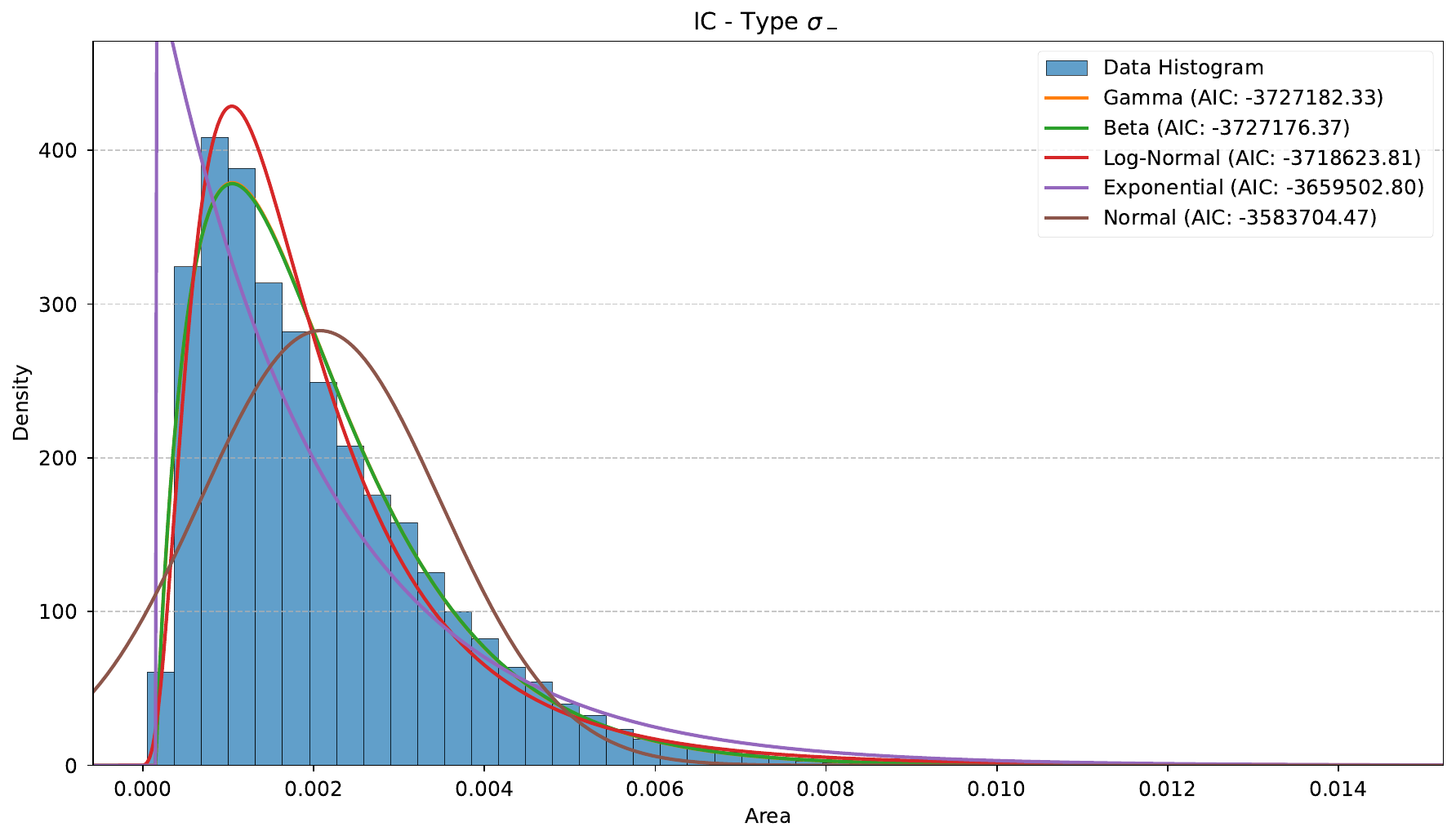}
        \caption{IC - Type $\sigma_-$}
        \label{fig:ic_sigma_minus_fits}
    \end{subfigure}

    \caption{
    Probability density functions of the areas for positive ($\sigma_+$) and negative ($\sigma_-$) terminal vortices. The square root of the most probable area for the EC case, $\sqrt{0.0125}$ corresponds to $1/10$ of the side length of the periodic domain. For the IC case, it corresponds to $3/100$ of the side length. The PDFs are fitted with various distribution functions.
    }
    \label{fig:all_histograms_fits_ec_ic}
\end{figure}

We will demonstrate that the statistical properties of terminal vortices can distinguish the two states (IC and EC).
We first observe distribution of the areas of the terminal vortices. 
The area of each terminal vortex is here defined as the number of pixels enclosed by the vortex divided by the total number of pixels $256 \times 256$.
The probability density functions (PDFs) of these areas are shown in the panels (a)--(d) of Figure~\ref{fig:all_histograms_fits_ec_ic}.
The PDF for the EC and IC states are calculated separately for $\sigma_+$ (counterclockwise rotating flows) and $\sigma_-$ (clockwise rotating flows). 
Within one state, the difference between the PDFs of $\sigma_+$ and $\sigma_-$ is small. 
In contrast, we observe a noticeable difference between the EC and IC states.
To quantify this difference, we fit the PDFs with various representative distributions of the random variables: normal, log-normal, gamma, beta, and exponential distributions. 
To evaluate goodness of fit, we adopt the Akaike Information Criterion (AIC). 
The AIC is minimized by the log-normal distribution for EC and the gamma distribution for IC,  regardless of $\sigma_+$ and $\sigma_-$, as indicated in Figure~\ref{fig:all_histograms_fits_ec_ic}. 
The log-normality of the EC state is consistent with the enstrophy-cascade picture in the sense of Kolmogorov 1962 \cite{Kolmogorov_1962}, provided that the area of the terminal vortex represents somehow the enstrophy flux (or dissipation rate) at the scale associated with the area. 
In contrast, the gamma distribution for the IC case may imply that the length scale associated with the area is a Gaussian random variable.

\begin{figure}[htbp]
    \centering

    \begin{subfigure}[b]{0.48\textwidth}
        \includegraphics[width=\textwidth]{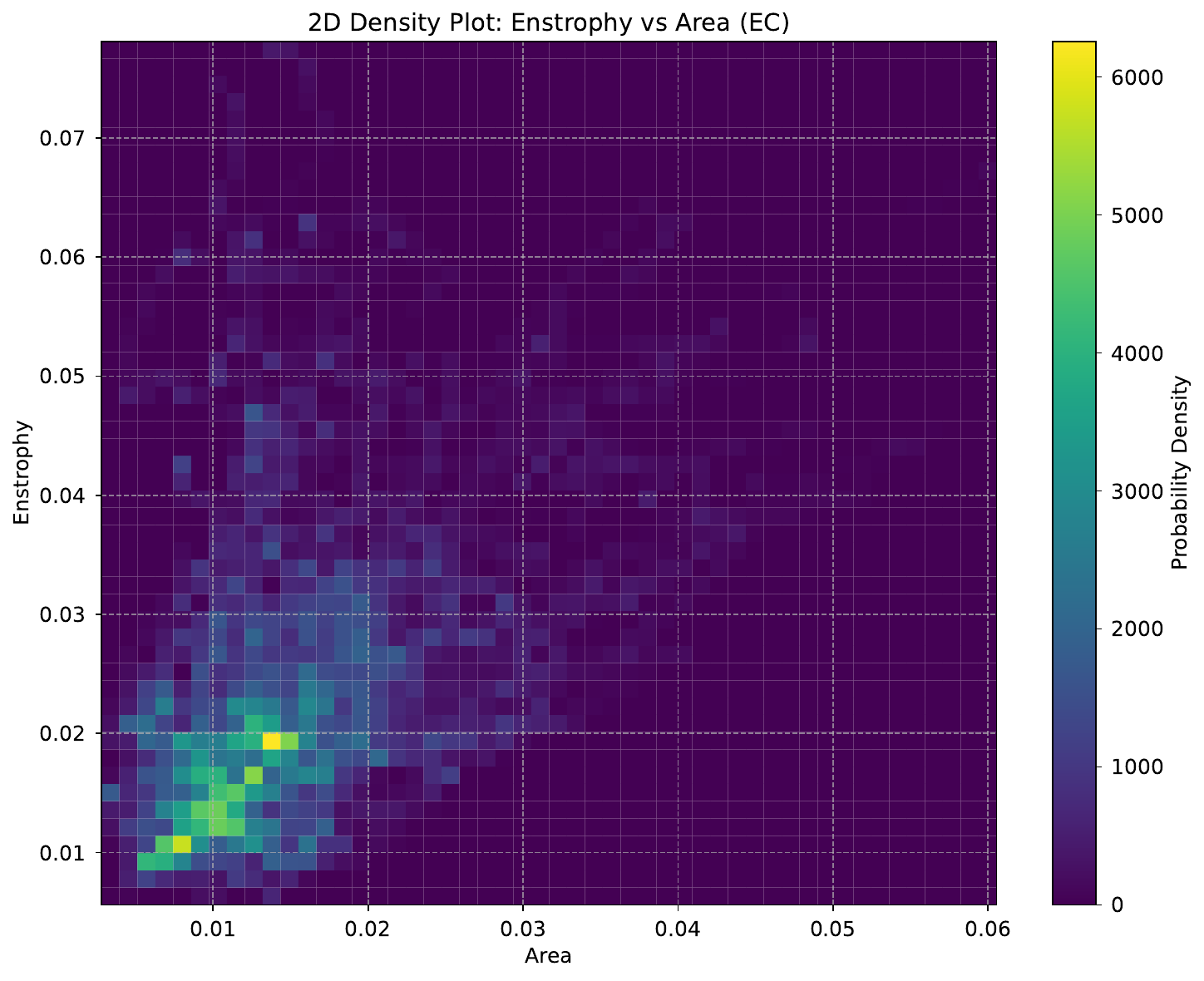}
        \caption{EC}
        \label{fig:2dhist_enstrophy_ec}
    \end{subfigure}
    \hfill
    \begin{subfigure}[b]{0.48\textwidth}
        \includegraphics[width=\textwidth]{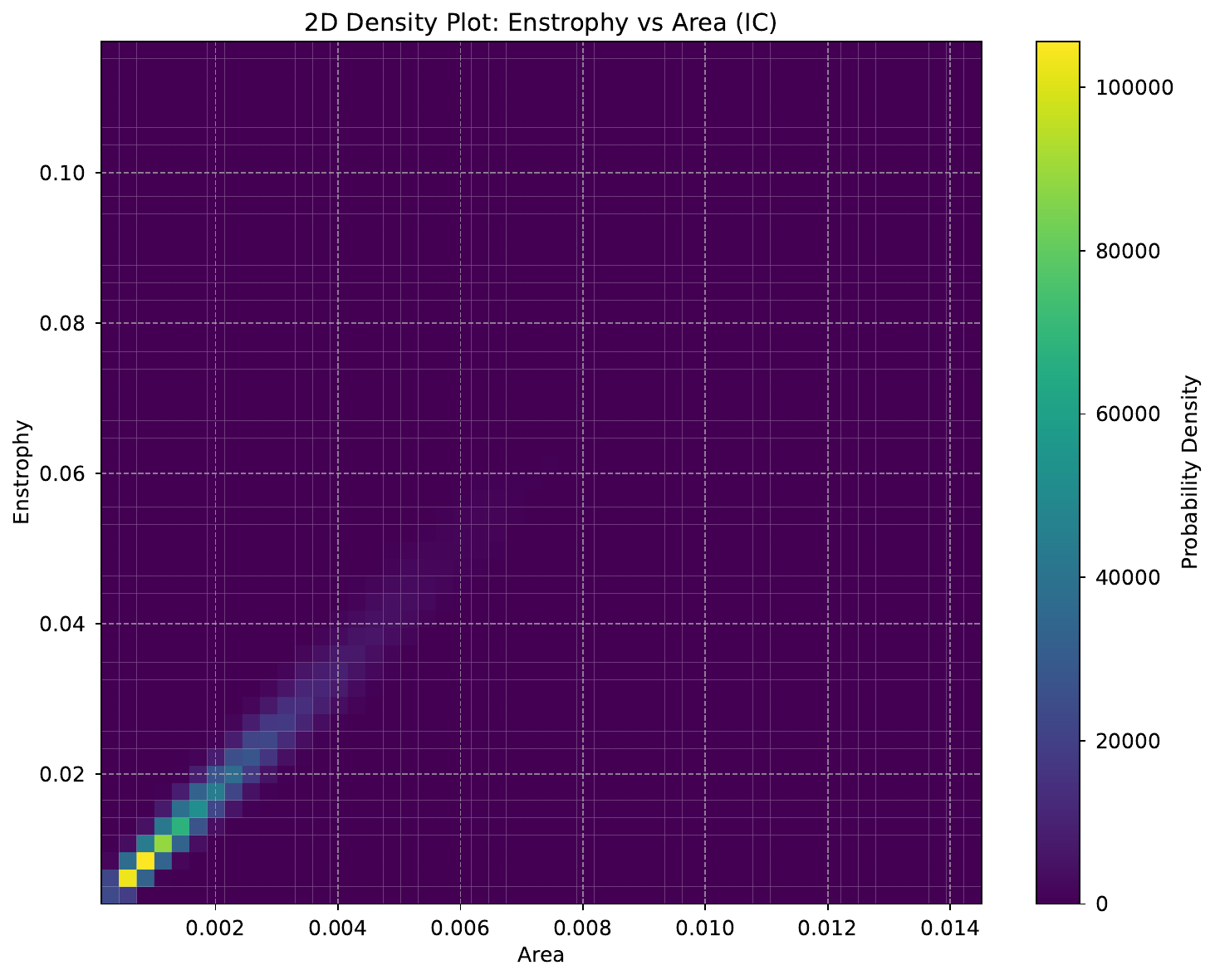}
        \caption{IC}
        \label{fig:2dhist_enstrophy_ic}
    \end{subfigure}
    \caption{Two dimensional joint PDFs of enstrophy versus vortex area for EC and IC data. Two-dimensional probability distribution function of the enstrophy and area of the terminal vortex for the EC and IC states.}
    \label{fig:2dhist_enstrophy}
\end{figure}

Motivated by these PDF of the area, we now study more directly relevant quantities to enstrophy or energy cascade through the terminal vortices.
We calculate the joint PDF of the area and the enstrophy on each terminal vortex.  
This enstrophy is computed by taking the sum of the squared vorticity on the coarse-grained grid points belonging to the same terminal vortex and dividing it by the total number of coarse-grained $256^2$ grid points. 
The joint PDFs of the EC and IC states shown in Figure~\ref{fig:2dhist_enstrophy} behave differently from each other.
The joint PDF of the EC state is broader than that of the IC state, indicating that the enstrophy and area of the EC state fluctuate more strongly than those of the IC state. 
For the IC case, the fluctuation is described by one variable (either the enstrophy or the area). 
Therefore, the joint PDF reveals a qualitative difference between the EC and IC states.
Now, the less fluctuating behavior of the joint PDF of the IC state can be understood by the following argument. 
For a given area $A$ of the terminal vortex, the typical enstrophy per area can be calculated as an integral of the enstrophy spectrum on the wavenumber from $A^{-1/2}$ up to the forcing wavenumber $k_f$, which is basically the largest wavenumber in the inertial range.
Let us now recall that the enstrophy spectrum, given by $k^2 E(k) \propto \epsilon^{2/3}k^{1/3}$ for the IC case, is an increasing function of $k$ in the inertial range. 
This implies that typical enstrophy is dominated by the contribution around $k_f$. 
Hence, no inertial range fluctuation is present in the joint PDF of the IC (Figure~\ref{fig:2dhist_enstrophy}), while inertial range fluctuations contribute in the joint PDF for the EC case.

\begin{figure}[htbp]
    \centering

    \begin{subfigure}[b]{0.48\textwidth}
        \includegraphics[width=\textwidth]{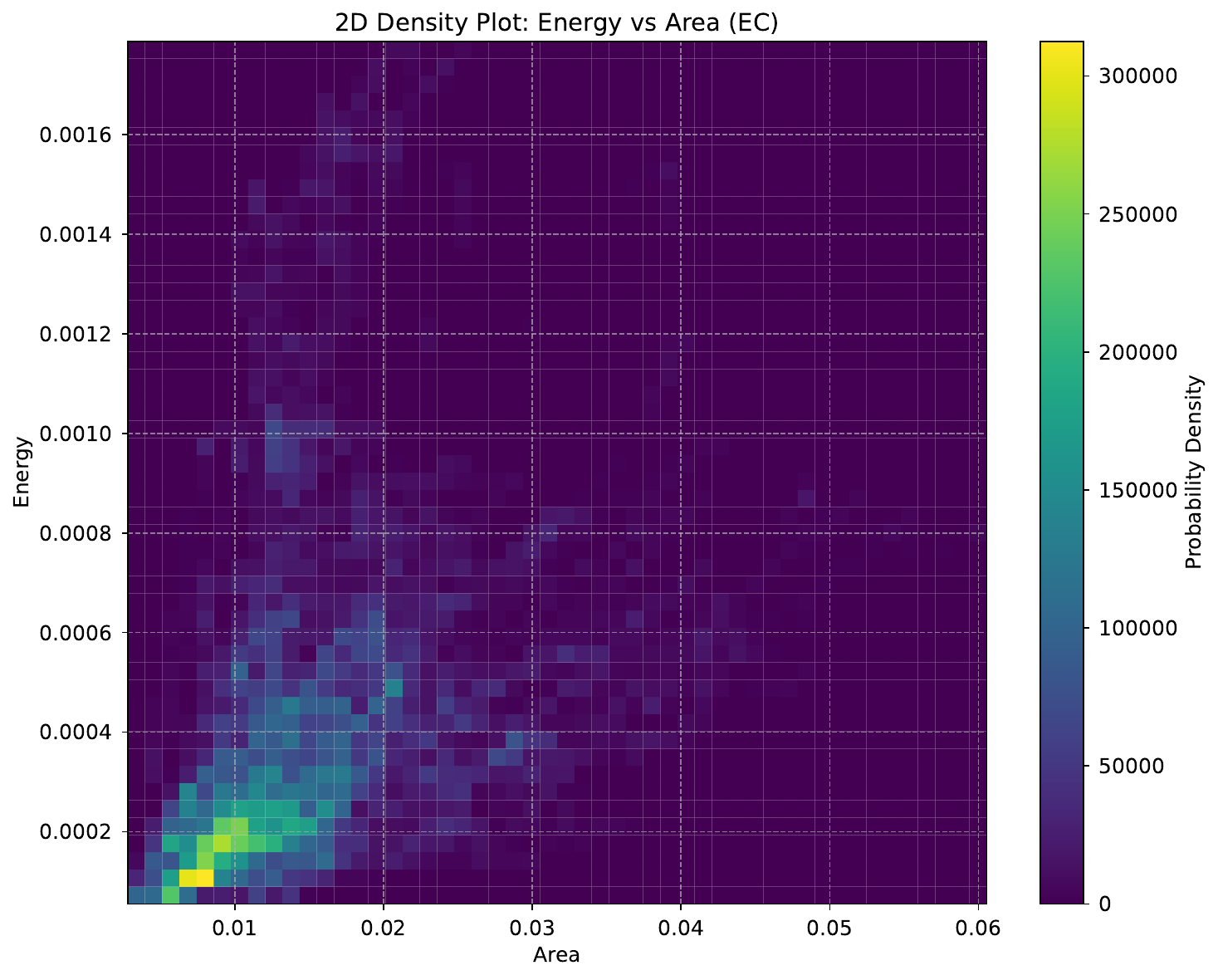}
        \caption{EC}
        \label{fig:2dhist_energy_ec}
    \end{subfigure}
    \hfill
    \begin{subfigure}[b]{0.48\textwidth}
        \includegraphics[width=\textwidth]{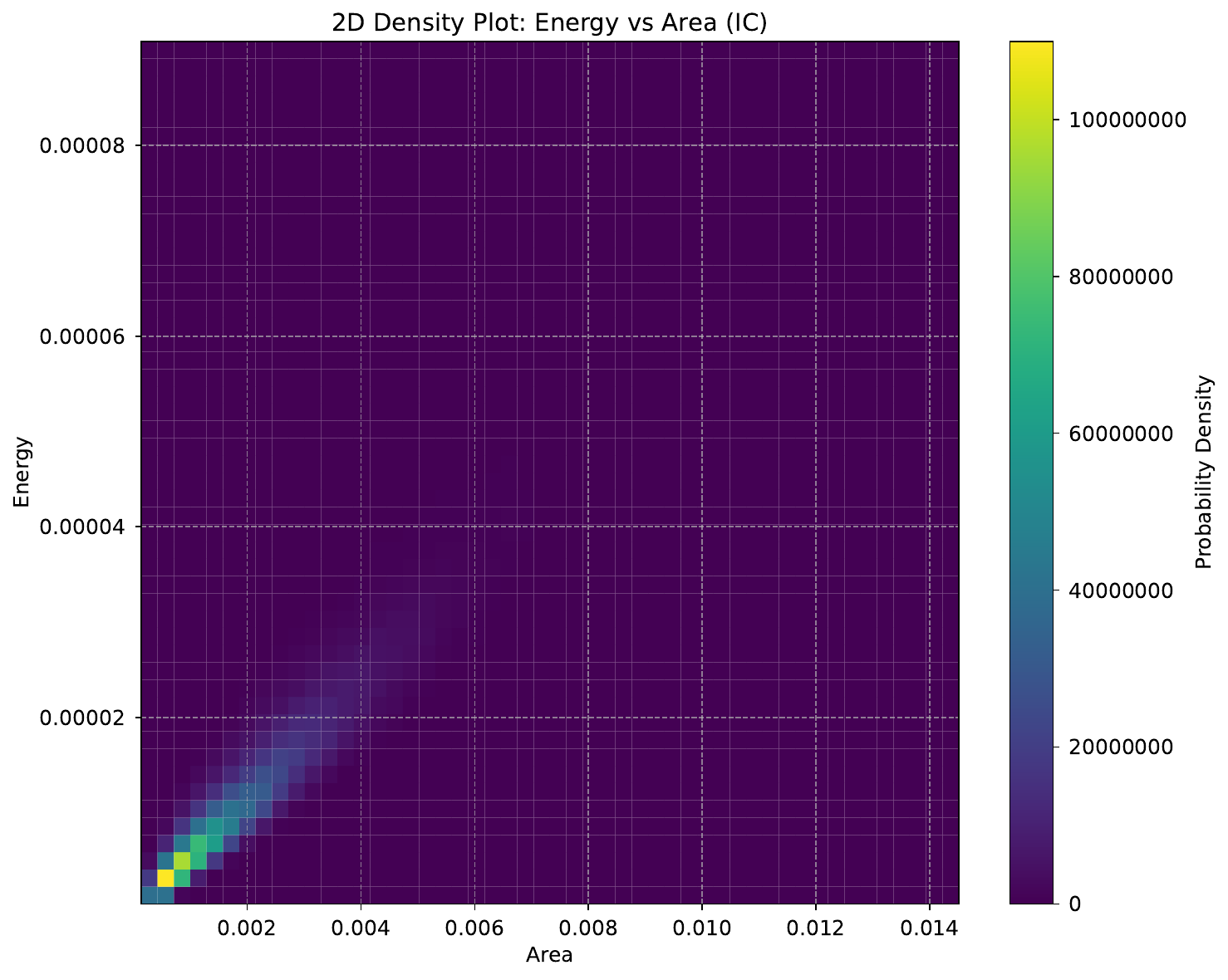}
        \caption{IC}
        \label{fig:2dhist_energy_ic}
    \end{subfigure}
    \caption{Two dimensional joint PDFs of energy versus vortex area for EC and IC data.}
    \label{fig:2dhist_energy}
\end{figure}

If we replace the enstrophy on the terminal vortex by the energy on it, we expect that the corresponding joint PDF for the IC state becomes two-dimensional.
However, as shown in Figure~\ref{fig:2dhist_energy}, the joint PDF of the energy and area of the IC state are one-dimensional, whereas those of the EC state are two-dimensional. 
The behavior is the same as that for the enstrophy-area joint PDFs. 
For the less fluctuating behavior of the energy--area joint fluctuations for the IC state, we do not have a simple explanation.
However, we just point out that the weaker fluctuations of the IC state are consistent with a common view of the absence of intermittency of the IC state~\cite{boffetta_two-dimensional_2012}.

In summary, we show here that the TFDA analysis is able to distinguish the two states of two-dimensional turbulence by utilizing the terminal vortices, more specifically looking at fluctuations of their areas, the enstrophy or the energy on them. 
The physical origin of the difference in the fluctuations of the EC and IC states will be reported elsewhere.

\section{Summary}\label{sec:6}
We have developed a topological classification theory for the orbit structures of structurally stable Hamiltonian flows on a torus.
Structurally stable Hamiltonian flows can have quasi-essential saddle separatrices such as Figure~\ref{fig:ess_saddle_connection}(b),
which differs from structurally stable Hamiltonian flows in the plane. 
Based on this classification theory, we show that the topological streamline structure of structurally stable Hamiltonian flows is converted into a partially cyclically ordered rooted tree (COT) when a torus has no physical boundaries. 
The conversion algorithm is summarized as follows. 
First, we reduce the Hamiltonian flow on the torus to a Hamiltonian flow on an annulus by cutting a handle of the torus along an essential periodic orbit. 
We then apply the conversion algorithm for structurally stable Hamiltonian flows on the annulus developed in~\cite{uda2019persistent_en}. 
The COT depends on this cut periodic orbit in general, but is uniquely determined by specifying the selection rule for the essential periodic orbit.
We note that we remove the loop from the Reeb graph of the Hamiltonian by cutting the torus, but two new nodes representing this periodic orbit are added to the COT.
However, the Reeb graph can be reproduced by identifying two nodes of the COT.

Finally, we apply the conversion algorithm to snapshots of freely decaying turbulence and enstrophy cascade turbulence over a planar domain with a doubly periodic boundary condition. 
Owing to this COT representation of the flows, coherent terminal vortex structures are extracted as regions bounded by quasi-inessential self-connected saddle separatrices. 
In addition, we examine the statistical properties of such terminal vortex structures and successfully distinguish between the energy cascade and the enstrophy cascade states in two-dimensional turbulence. 
The physical mechanism of such statistical states in terms of the interactions of the terminal vortices will be reported soon.

\appendix
\section{Procedure (P) choosing an essential periodic orbit}\label{sec:appA}
\begin{figure}[t]
\begin{center}
\includegraphics[bb=0 0 740 287,scale=0.475]{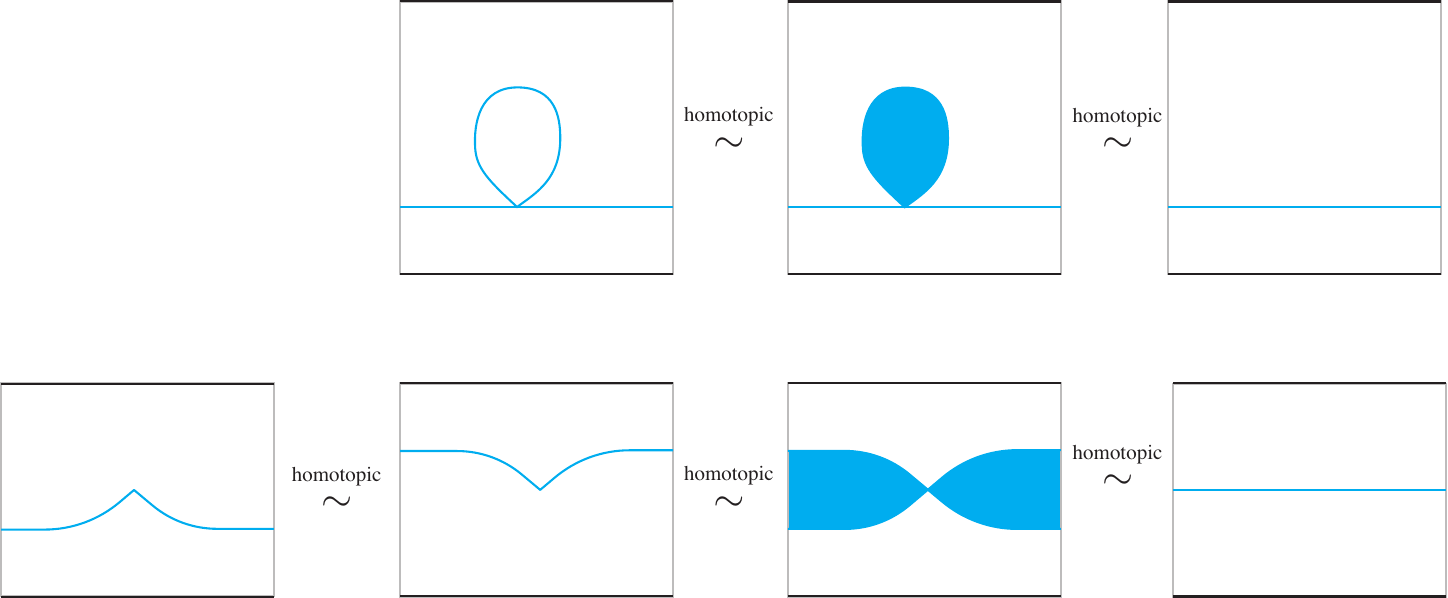}
\end{center}
\caption{First homology class of saddle connections for the local orbit structures  (a) $a_{-\dot{\pm}}$ and (b) $\alpha_{+\dot{-}}$.
The symbol $\sim$ means that the subgroups in $H_1(\mathbb{T}^2)$ induced by 
the inclusions coincide with each other and are generated by one generator.}
\label{fig:homology_class}
\end{figure}
Let us recall some mathematical terms and known facts.
We consider the local orbit structure of $a_{-\dot{\pm}}$ shown in the top of Figure~\ref{fig:homology_class}(a).
When we fill the disk bounded by the inessential saddle connection $C$, the saddle connection $C$ and its filling $F_C$ are homologous and belong to the same homology class. 
More precisely, the first homologies induced by inclusions $C \to S$ and $F_C \to S$ are isomorphic and are generated by one element respectively. 
Furthermore, the generators are the same homology class except for the sign as one essential simple closed curve on the right panel of Figure~\ref{fig:homology_class}(a).
Similarly, as shown in Figure~\ref{fig:homology_class}(b), the two essential simple closed curves in the local orbit structure of $\alpha_{-\dot{+}}$ are homologous to the one obtained by filling a contractible disk region between the saddle connections. 
Hence, they all belong to the same homology class
as the essential simple closed curve in the rightmost ones of Figure~\ref{fig:homology_class}(a,b).
On the other hand, every simple closed curve $\gamma$ in the torus $\mathbb{T}^2:= \R^2 / \Z^2$ can be winding $k$ times horizontally and $l$ times vertically as an undirected circuit. 
More precisely, a simple closed curve $\gamma$ defines the projection class $[(k,l)] \in (\{ 0 \} \times \Z_{> 0}) \sqcup (\Z_{> 0} \times \Z) \cong (\Z^2 - \{ 0 \})/\pm \cong (H_1(\mathbb{T}^2; \Z) - \{ 0 \}))/\pm$ of the first homology class.
In particular, $a_{-\dot{\pm}}$ and $\alpha_{-\dot{+}}$ in Figure~\ref{fig:homology_class}(a,b)
define the first homology class corresponding to $[(1,0)]$.
\begin{figure}[t]
\begin{center}
\includegraphics[bb= 0 0 209 217, scale=0.65]{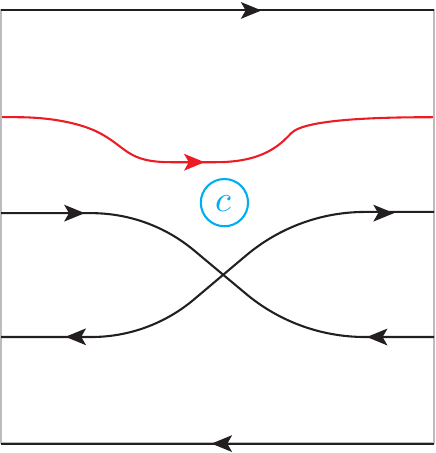}
\end{center}
\caption{An essential periodic orbit (a level curve of the Hamiltonian) near the saddle connection of the local orbit structure $\alpha_{-\dot{+}}$.}
\label{fig:cutting_loop_01}
\end{figure}

\begin{figure}[t]
\begin{center}
\includegraphics[bb=0 0 477 107, scale=0.65]{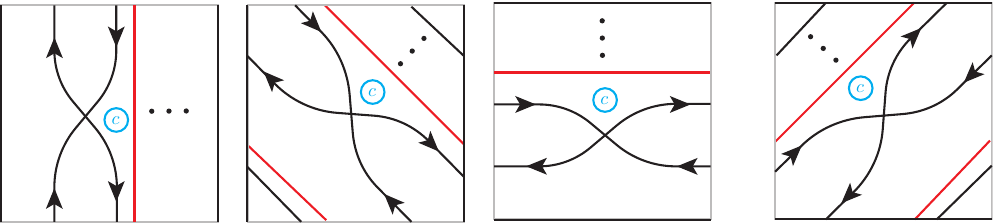}
\end{center}
\caption{First homology class of saddle connection}
\label{fig:cutting_loop_02}
\end{figure}

Using these terms, the procedure (P) for defining an essential periodic orbit is now precisely described as follows:
The projection class of the first homology class for the essential simple closed curve $\gamma$ in the saddle connection of the pattern $\alpha_{-\dot{+}}$ in Figure~\ref{fig:cutting_loop_01} is $[(k,l)]=[(1,0)]$. 
Therefore, the essential periodic orbit is chosen by translating in the direction perpendicular to the $[(k,l)]$ direction within the level set near $\gamma$ shown as the red orbit in Figure~\ref{fig:cutting_loop_01}.
Here we note that there are two choices of essential simple closed curves as shown in the two left panels of Figure~\ref{fig:homology_class}.
On the other hand, regardless of which of these two is chosen, the same projection class of the first homology class is determined.
Hence, for the saddle connection $\gamma$ with $[(k,l)] \in (\{ 0 \} \times \Z_{> 0}) \sqcup (\Z_{> 0} \times \Z)$, 
we define an essential periodic orbit as follows.
\begin{enumerate}
\item For $k=0$,  since the essential simple closed curve does not wind horizontally,  and the simple closed curve $\gamma$ has no self-intersection, we have $[(k,l)] = [(0,1)]$ as shown in Figure~\ref{fig:cutting_loop_02}(a). Then, we choose a level set near $\gamma$ shifted along the $x$-axis as an essential periodic orbit.
\item For $k>0$, as shown in Figure~\ref{fig:cutting_loop_02}(b,c,d), we choose a level set near $\gamma$ shifted along the $y$-axis as essential periodic orbit.
\end{enumerate}

\begin{remark}
Note that the essential periodic orbit extracted by this procedure is not continuous for the choice of projection class for the first homology class.
\end{remark}

\noindent
{\bf Acknowledgments}
We thank the members of Technology Research \& Innovation, BIPROGY Inc. for useful discussions and development of the Python library, $\mathbb{T}^2$-\texttt{psiclone}, used in this research.
The research is partially supported by JST MIRAI JPMJMI22G1.

\bibliographystyle{abbrv}
\bibliography{KSY}

\end{document}